\documentclass[prepring, 11pt, authoryear, 3p]{elsarticle}
\usepackage[american]{babel}

\usepackage[titletoc,toc,title]{appendix}
\usepackage{hyperref}

\makeatletter
\def\ps@pprintTitle{%
 \let\@oddhead\@empty
 \let\@evenhead\@empty
 \def\@oddfoot{}%
 \let\@evenfoot\@oddfoot}
\makeatother

\usepackage{multicol}

\usepackage{pst-plot}
\usepackage{subfig}

\usepackage{makecell}

\usepackage{algorithm2e}

\usepackage{url}

\usepackage{amsmath}
\usepackage{amssymb}
\usepackage{amsthm}
\usepackage{amsfonts}

\usepackage{multirow}

\usepackage{graphicx}
\DeclareGraphicsExtensions{.eps,.pdf,.jpg,.png}

\usepackage{multicol}

\renewcommand\vec{\mathbf}

\newtheorem{definition}{Definition}
\newtheorem{theorem}{Theorem}
\newtheorem{lemma}{Lemma}

\newcommand*\nd{{\mkern 1mu\cdot\mkern 1mu}}

\usepackage{array}

\usepackage{xcolor}
\usepackage{colortbl}
\newcommand{\cgf}{\cellcolor{lightgray!100}}
\newcommand{\cgs}{\cellcolor{lightgray!60}}
\newcommand{\cgt}{\cellcolor{lightgray!20}}

\begin{document}

\begin{frontmatter}

\title{An Extension of the Non-Inferior Set Estimation Algorithm\\ for Many Objectives}

\ifx
\author[add1]{Author1}
\ead{email1@domain.ca}
\author[add2,add3]{Author2}
\ead{email2@domain.ca}
\author[add2,add3]{Author3}
\ead{email3@domain.ca}

\address[add1]{Department of Redundancy Department}
\address[add2]{Centre for Study of Things}
\address[add3]{Department of Interests}
\fi

\cortext[cor1]{Corresponding author}
\author{Marcos M. Raimundo\corref{cor1}}
\ead{marcosmrai@gmail.com}
\author{Paulo A. V. Ferreira}
\ead{valente@dt.fee.unicamp.br}
\author{Fernando J. Von Zuben}
\ead{vonzuben@dca.fee.unicamp.br}
\address{University of Campinas - UNICAMP, School of Electrical and Computer Engineering, Av. Albert Einstein -- 400, 13083-852, Campinas, S\~ao Paulo, Brazil}

\begin{abstract}

This work proposes a novel multi-objective optimization approach that globally finds a representative non-inferior set of solutions, also known as Pareto-optimal solutions, by automatically formulating and solving a sequence of weighted sum method scalarization problems. The approach is called MONISE (Many-Objective NISE) because it represents an extension of the well-known non-inferior set estimation (NISE) algorithm, which was originally conceived to deal with two-dimensional objective spaces. 
The proposal is endowed with the following characteristics: (1) uses a mixed-integer linear programming formulation to operate in two or more dimensions, thus properly supporting many (i.e., three or more) objectives; (2) relies on an external algorithm to solve the weighted sum method scalarization problem to optimality; and (3) creates a faithful representation of the Pareto frontier in the case of convex problems, and a useful approximation of it in the non-convex case. Moreover, when dealing specifically with two objectives, some additional properties are portrayed for the estimated non-inferior set. Experimental results validate the proposal and indicate that MONISE is competitive, in convex and non-convex (combinatorial) problems, both in terms of computational cost and the overall quality of the non-inferior set, measured by the acquired hypervolume.

\end{abstract}


\end{frontmatter}

\section{Introduction}

Many practical applications are better modeled as optimization problem, characterized by the existence of multiple conflicting objectives. A classical and usual example is the compromise between maximizing consumer satisfaction and minimizing service cost. Indeed, dealing with conflicting objectives is omnipresent in our lives, and a significant portion of these multi-objective problems admits a proper mathematical formulation, so that we may resort to computational resources to obtain Pareto-optimal solutions, also called non-inferior solutions \citep{Miettinen1999}. 

Obviously, the main challenge of multi-objective optimization is the need of simultaneously dealing with conflicting objectives. Given the multidimensional nature of the objective function, two solutions $\vec{y}$ and $\overline{\vec{y}}$ only establish a dominance relation with each other when all the components of the solution $\vec{y}$ are equally or better satisfied in comparison to what happens in the case of solution $\overline{\vec{y}}$, with $\vec{y}$ being strictly better in at least one objective. We are going to properly define the most relevant multi-objective concepts in the next section.

Most solution techniques to multi-objective optimization have been conceived to deal with problems characterized by two or at most three conflicting objectives. The extension to more objectives is not straightforward in some cases and is even not possible in other situations. Besides, with the increase in the number of objectives, scalability issues arise, together with a dramatic reduction in the relevance of the concept of dominance \citep{Kukkonen2007, Ishibuchi2009}, thus imposing amazing challenges to algorithms based mainly on dominance relations \citep{Deb2002,Zitzler2001}. To overcome this issue, many-objective population-based algorithms have been proposed \citep{Deb2013, Ishibuchi2009}. They generally rely on scalarization-based approaches such as the weighted sum method \citep{Marler2009}, reference points \citep{Das1998} and box-constrained models \citep{Caballero2004} to build algorithms capable of producing a consistent approximation of the Pareto frontier.

In contrast to the early-conceived multi-objective heuristically-based algorithms, scalarization is the cornerstone of many deterministic \textit{a posteriori} algorithms for multi-objective optimization, and the main proposals are as follows: \cite{Das1998} proposed a method in which well-spaced points are calculated inside the hyperplane supported by the individual minima.  Based on that, collinear solutions (forced using equality constraints) outlined by these points and the normal vector to this hyperplane are searched. \cite{Messac2003} proposed a similar process using inequality constraints, making the solution as collinear as possible, according to the optimization process. Further adjustments of this method were made by \cite{Messac2004} and \cite{Sanchis2007}. Other methods also used some initial pieces of information to calculate a \textit{set of parameters} for scalarization, thus finding all the aimed representations in parallel \citep{Snyder1997,Burachik2013,Khorram2014}. On the other hand, adaptive methods resort to already known information about the Pareto frontier to iteratively find new efficient solutions: \cite{Ryu2009a} iteratively determined the solution that is the farthest from its neighbors, creating a second order approximation using its neighbors and optimizing this approximation inside a thrust region radius; \cite{Ozlen2009} and \cite{Ozlen2014} proposed a recursive algorithm that, setting a superior limit to the $k$-th objective (using information provided by already found solutions), recursively applies the same method to solve the reduced problem composed of $k-1$ objectives; \cite{Sylva2004} used a linear-integer programming formulation to exclude regions dominated by already found solutions; \cite{Masin2008a} used a linear-integer programming formulation to find efficient solutions maximizing the distance to other already found solutions; \cite{Eichfelder2009, Eichfelder2009a} resorted to the already found efficient solutions to create a first-order approximation used to determine new parameters for the Pascoaletti-Serafini scalarization; \cite{Kim2005} proposed a method that initially creates a rough representation of the Pareto frontier, further prospecting poorly explored regions by finding solutions that are collinear with the segment determined by the Nadir point and an ``expected'' solution estimated at a poorly explored region. An attempt to split the search region into boxes, excluding regions dominated by already found efficient solutions as well as regions that would dominate them, characterized the work of \cite{Kirlik2014c}, in which rough boxes are centered on already found solutions, and also the work of \cite{Ceyhan2019}, in which non-dominated reference points and regions around already found solutions are recursively created.

In the class of algorithms that use outer and/or inner representations to approximate the Pareto frontier, the Benson's algorithm \citep{Benson1998} constructs a representation of linear programming problems by constructing a simplex in the objective space that is refined by recursively finding vertices from the simplex that do not have a corresponding solution in the feasible set, and adding hyperplanes that separate the current vertex, thus refining the simplex.

The Benson's algorithm \citep{Benson1998} was further improved to $\epsilon$-efficiency \citep{Shao2008} by not separating vertices whose distance to the closest feasible solution is lower than $\epsilon$. This work also simplifies some steps of the algorithm by reducing the number of linear-programming steps; \cite{Ehrgott2011} extension approached $\epsilon$-approximated multi-objective convex optimization problems using the first order information to create the hyperplanes. Instead of using hyperplanes to construct a simplex, \cite{Nobakhtian2017} extension approaches nonconvex problems by using cones to construct an inner approximation of dominated cones and an outer approximation cone that would dominate already found efficient solutions.

The main problem for those methods resides on the type of scalarization they employ; which necessarily increase the complexity of the problem by recursively imposing additional constraints \citep{Messac2003, Ryu2009a, Eichfelder2009a, Eichfelder2009, Kim2005}, which in turn increases the complexity of solving each scalarization step. One way to bypass this problem relies on the weighted sum method scalarization. It consists only in making a linear combination of the objectives, thus not adding any constraint, and hence keeping the same complexity of the original problem. One of the striking algorithms that rely on the weighted sum method is the Non-Inferior Set Estimation (NISE) method \citep{Cohon1979}. NISE iteratively finds a solution using a scalarization based on the weighted sum method, ensuring a low computational cost and a simple implementation. It also sets up an inner and an outer representation of the frontier by connecting neighboring solutions with straight lines to build the inner representation and uses the hyperplanes formed by the optimized solutions to build the outer representation.

However, the recursive application of the weighted sum method by NISE will fail when dealing with problems having three or more objectives. Some extensions for higher dimensions, called polyhedral approximation approaches, instead of connecting neighboring solutions to build the inner representation, enumerate all facets of the convex hull and use the facets of the hull as weights to perform the optimization of the next weighted sum method \citep{Solanki1993, Klamroth2003, Craft2006, Rennen2009, Bokrantz2013}. Since a facet can have a non-positive component, \cite{Solanki1993} calculated the importance of the next facet by finding the larger distance to a point constrained to the outer approximation, which consists of all hyperplanes of the optimized solutions by adding a constraint not to allow the optimization to go further than the worst possible objective; \cite{Klamroth2003}, instead of using the facets of the convex hull, used its vertices as references to build an inner approximation with maximal combination of the vertices of each facet or a solution closer to any vertex in the outer approximation (thus not using the weighted sum method); \cite{Craft2006} assessed the importance of every facet by calculating the distance from the facet to the lower distal point, which is the point provided by the combination of the hyperplanes of the solutions from the vertices of the facet; the non-positivity of the facets is dealt with a heuristic that finds the most distinct combination of the weights from the vertices; \cite{Rennen2009} used the same approach of \cite{Solanki1993} to calculate the importance of every facet but introduced the concept of dummy-points which act as additional points at the convex hull, besides the points found by the optimization; to deal with mixed integer problems, \cite{Ozpeynirci2010a} resorted to an approach similar to the one proposed by \cite{Rennen2009}. Those points prevent the facets from having non-positive components. \cite{Bokrantz2013}, instead of finding the convex hull of the inner approximation, built the convex hull of the outer approximation, calculated the importance of every vertex from the outer approximation, and found the corresponding weight for that vertex. Also, \cite{Bokrantz2013} only calculated the importance of a vertex if the previous importance (less constrained) was higher than the importance of already known vertices.

Supported by theoretical evidence and focusing on the interplay of the weighted sum method and the NISE algorithm, we propose here an extension of NISE for more than two objectives, called MONISE. Consistent with the original proposal, our method is composed of an \textit{a posteriori} adaptive method which resorts to a set of solutions, together with their corresponding weighting vectors, to recursively calculate the next weighting vector to be used by the solver to achieve a new efficient solution. MONISE involves a standard mixed-integer linear programming formulation, for which commercial and open-source efficient solvers are available. The main distinctive aspect of our approach resides in the mixed-integer linear programming formulation of a general multi-objective optimization problem. Basically, it simplifies the search for the best facet in the convex hull by cleverly searching only among the most promising facets/outer vertices of the inner/outer approximations. Therefore, it avoids enumerating all facets of the whole convex hull and also avoids unnecessary steps of linear programming done in \cite{Solanki1993} and \cite{Rennen2009} to calculate the importance of every facet. \cite{Bokrantz2013} cleverly noted the existence of unnecessary calculations, but their approach to avoid them differs from the one proposed in this paper, since it only avoids recalculations of already enumerated facets, while our approach only explores and calculates the importance of the most promising facets.

Before presenting and experimentally validating MONISE, we perform an investigation of the weighted sum method and of the NISE algorithm. After presenting our proposal, we will demonstrate that MONISE produces consistent behaviour when two or more conflicting objectives are considered, and represents a scalable proposal when the number of objectives significantly increases.

This work is organized as follows: Section \ref{sec:def} formally describes the main concepts of multi-objective optimization; Section \ref{ssec:weighted_problem} deals with the weighted sum method and its main properties; Section \ref{ssec:nise} delineates the Non-Inferior Set Estimation adaptive algorithm, and Section \ref{ssec:moo:nise_disc} evidences NISE characteristics and properties;  Section \ref{sec:monise} presents our extension of NISE for many objectives (denoted MONISE), its characteristics, properties and main distinctions when compared with the original NISE; Section \ref{sec:experiments} is devoted to the description of some experiments toward validating MONISE and assessing its potential for multi-objective optimization, specially for a large number of objectives; Section \ref{sec:conclusion} summarizes the work with an analytical view of the findings, also including future perspectives of the research. Looking for a more fluent reading, all the proofs of pertinent theorems are left as Appendices A to D.

\section{Conceptual aspects of multi-objective optimization} \label{sec:def}

Let us firstly define a multi-objective problem.

\begin{definition}\label{def:moo:opt_mo}
A multi-objective problem is defined as follows \citep{Marler2004}:

\begin{equation*}
\begin{aligned}
& \underset{\vec{x}}{\text{minimize}} 
& &  \vec{f}(\vec{x}) \equiv \{f_1(\vec{x}), f_2(\vec{x}),\ldots, f_m(\vec{x})\}\\
& \text{subject to}
& & \vec{x} \in \Omega,\ \Omega \subset \mathbb{R}^{n},\\
& & & \vec{f}(\cdot): \Omega \rightarrow \Psi,\ \Psi \subset \mathbb{R}^m.
\end{aligned}
\end{equation*}
\end{definition}

\begin{figure}[ht!]
\centering
\includegraphics[width=0.75\textwidth]{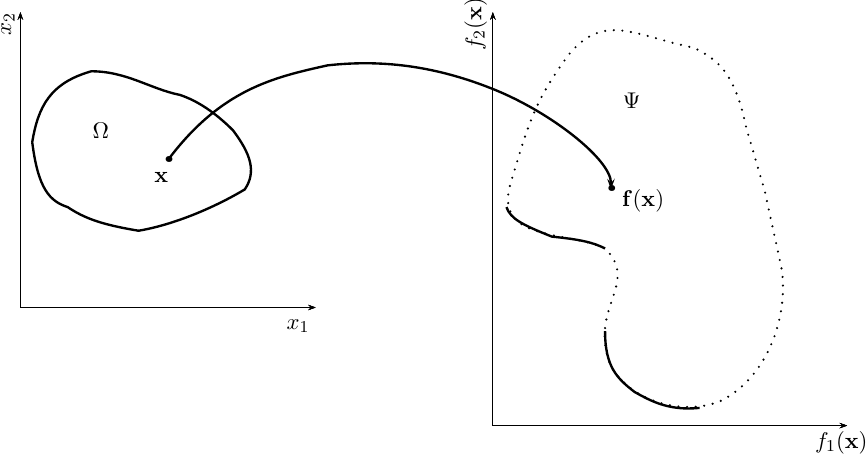}
 \caption{Representation of the decision space (on the left) and of the objective space (on the right), taking two decision variables and two objectives.}
 \label{fig:moo:espacos_mo}
\end{figure}

In Definition \ref{def:moo:opt_mo}, the set $\Omega\subset \mathbb{R}^n$ is known as the decision space and the set $\Psi \subset \mathbb{R}^m$ is known as the objective space. Figure \ref{fig:moo:espacos_mo} represents a relationship between those two spaces (restricted to two dimensions for visualization purposes). Each point at the decision space has a correspondent point at the objective space, obtained by evaluating each objective function. On the objective space, the two bold lines correspond to the Pareto frontier, which is the set of all efficient or non-inferior solutions.

In the sequence, based on the formalism provide by \citet{Marler2004}, we present some basic definitions to contextualize the multi-objective optimization problem. Without loss of generality, the objectives are associated with minimization problems.

\begin{definition}\label{def:moo:efi} \textbf{Efficiency/Pareto-optimality}:
A solution $\vec{x}^*\in \Omega$ is efficient (Pareto-optimal) if there is no other solution $\vec{x} \in \Omega$ such that $f_i(\vec{x}) \leq f_i(\vec{x}^*)$, $\forall i \in \{1,2,\ldots,m\}$ and $f_i(\vec{x}) < f_i(\vec{x}^*)$ for at least one $i \in \{1,2,\ldots,m\}$.
\end{definition}

\begin{definition}\label{def:moo:efi_front} \textbf{Efficient frontier/Pareto frontier}:
The efficient frontier $\Psi^{*}$ (Pareto frontier) is the set of all efficient solutions. When considered the problem in Definition \ref{def:moo:opt_mo}, the efficient frontier $\Psi^{*}$ is formed by efficient objective vectors $\vec{f}(\vec{x}^*) \in \Psi^{*}$ which have corresponding feasible solutions $\vec{x}^* \in \Omega$. Also, $\Omega^{*}$ is the set of feasible solutions whose objective vectors belong to the efficient frontier: $\vec{x}^* \in \Omega^{*} \Leftrightarrow \vec{f}(\vec{x}^*) \in \Psi^{*}$.
\end{definition}

The following definitions are necessary to support the proposition of some of the adaptive and scalarization methods. The ``$k$-th definitions'' are intended to refer to single objective solutions.

\begin{definition}\label{def:moo:indiv_min_value} \textbf{$k$-th individual minimum value}:
When only the $k$-th component of the objective function vector is optimized, a solution $\vec{x}^{*(k)}$ is obtained. The $k$-th individual minimum value $l^{(k)}$ corresponds to the minimum value of the optimization $(l^{(k)} = f_k(\vec{x}^{*(k)}))$.

\begin{equation*}
\begin{aligned}
& \underset{\vec{x}}{\text{minimize}}
& & f_k(\vec{x})\\
& \text{subject to}
& & \vec{x} \in \Omega,\ \Omega \subset \mathbb{R}^{n},\\
& & & \vec{f}(\cdot): \Omega \rightarrow \Psi,\ \Psi \subset \mathbb{R}^m.
\end{aligned}
\end{equation*}
\end{definition}

\begin{definition}\label{def:moo:indiv_min_sol} \textbf{$k$-th individual minimum solution}:
An individual minimum solution $\vec{l}^{*(k)}$ is an efficient solution characterized by having its $k$-th component equal to the $k$-th individual minimum value $l^{(k)}$.
\end{definition}

\begin{definition}\label{def:moo:utop_sol} \textbf{Utopian solution}:
A utopian solution $\vec{z}^{utopian}$ is a vector in the objective space characterized by having its $k$-th component $z^{utopian}_k$ given by the $k$-th individual minimum value $l^{(k)}$ (see Definition \ref{def:moo:indiv_min_value}), and this is valid for all $k \in \{1,2,\ldots,m\}$:
\begin{equation*}
\vec{z}^{utopian}=\{l^{(1)};\ldots;l^{(m)}\}.
\end{equation*} 
It is important to notice that, in case of conflicting objectives, the utopian solution is not an attainable solution.
\end{definition}

\section{The weighted sum method and the Non-inferior Set Estimation} \label{sec:nise}

\subsection{The weighted sum method} \label{ssec:weighted_problem}

The weighted sum method consists in optimizing a convex combination of the objectives, with each component of the weighting vector representing a relative importance of the corresponding objective. With this scalarization, the designer aims at expressing his/her preferences, given the objectives \citep{Cohon1979}. Additionally, as will be done in this work, the weighting vector may be automatically determined by a recursive process, with the purpose of exploring particular regions of the Pareto frontier.

\begin{definition}\label{def:moo:weight_prob}
The definition of the weighted sum method is given by:

\begin{equation*}
\begin{aligned}
& \underset{x}{\text{minimize}}
& & \vec{w}^\top \vec{f}(\vec{x})\\
& \text{subject to}
& & \vec{x} \in \Omega,\ \Omega \subset \mathbb{R}^n,\\
& & & \vec{f}(\cdot): \Omega \rightarrow \Psi, \Psi \subset \mathbb{R}^m.
\end{aligned}
\end{equation*}
where $\sum_{i=1}^m w_i = 1,\ \vec{w} \in \mathbb{R}^m$ and $w_i \geq 0\ \forall i \in \{1,2,\ldots,m\}$.
\end{definition}

\begin{figure}[ht!]
\centering
\includegraphics[width=0.5\textwidth]{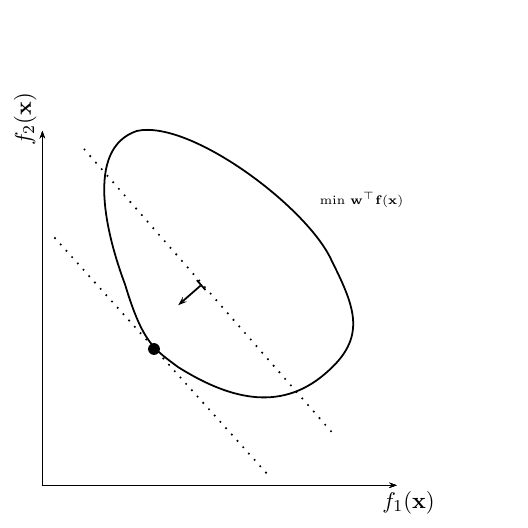}
\caption{Representation of the solution produced by the weighted sum method.}
 \label{fig:moo:weight_prob}
\end{figure}

In Figure \ref{fig:moo:weight_prob}, the weighting vector $\vec{w}$ defines the slope of the line that guides the optimization process, reaching a tangent point in the objective space.

Some properties of this scalarization are relevant. In the general case, without assuming any particularity of the objective space $\Psi$, an optimal solution for the weighted sum method results in an efficient solution \citep{Geoffrion1968, Miettinen1999}, and all efficient solutions are attainable by the weighted sum method if the problem is convex \citep{Miettinen1999}. In the non-convex case, in general, any efficient solution $\vec{f}(\overline{\vec{x}})$ which does not admit supporting hyperplanes to $\Psi$ passing through them are not attainable by the weighted sum method \citep{Koski1985, Das1997}, which means that there is no weight vector $\vec{w}$ capable of conducting the weighted sum method to find an efficient solution $\overline{\vec{x}}$ with objective vector $\vec{f}(\overline{\vec{x}})$.

\subsection{NISE - Non-inferior Set Estimation} \label{ssec:nise}

The NISE (Non-inferior Set Estimation) method \citep{Cohon1979} is an iterative approach that uses the weighted sum method to automatically create, at the same time, an inner and an outer approximation of the Pareto frontier using a linear approximation. At every iteration, based on the already calculated efficient solutions, it is traced a segment between each neighboring pair of solutions, determining new weighting vectors. This procedure finds an accurate and fast approximation for problems with two objectives \citep{Romero2003}.

Two neighboring efficient solutions (called neighborhood) are used to determine a new efficient solution employing the weighted sum method. Specifically: the initialization should generate the first two solutions (Section \ref{sssec:moo:nise:ini}); at each iteration, the next neighborhood to be explored should be determine (Section \ref{sssec:moo:nise:choose}), thus obtaining the parameters for the weighted sum method (Section \ref{sssec:moo:nise:w_inf}), followed by a new solution, and new neighborhoods (Section \ref{sssec:moo:nise:branch}). The stopping criterion is defined to ensure the quality threshold of the approximation (Section \ref{sssec:moo:nise:stop_crit}).

\subsubsection{Initialization} \label{sssec:moo:nise:ini}

The initialization consists in finding the first two solutions $\vec{r}^{*1}=[f_1(\vec{x}^{*1}),$ $ \ldots, f_m(\vec{x}^{*1})]$ and $\vec{r}^{*2}=[f_1(\vec{x}^{*2}), \ldots, f_m(\vec{x}^{*2})]$ which are individual minimum solutions (Definition \ref{def:moo:indiv_min_sol}) for objectives $1$ and $2$, respectively. The weighting vectors $\vec{w}^{1}$ and $\vec{w}^{2}$ have null elements except for the element corresponding to the objective being optimized, assumed to be equal to 1. Finally, it is possible to define the first neighborhood $\mathcal{N}^1 = \{(1,2)\}$ which will be used to find the subsequent solutions.

\subsubsection{Neighborhood choice} \label{sssec:moo:nise:choose}

Considering a set of neighborhoods $\mathcal{N}^k$, the neighborhood to be explored at the $k$-th iteration is the one that has the maximum error $\mu = \max\ \mu^{i,j}, \forall (i,j) \in \mathcal{N}$. To find that error, for every neighborhood $(i,j) \in \mathcal{N}$ it is possible to calculate the largest distance $\mu^{i,j}$ between the normal vector $\vec{w}$ (calculated as described in Section \ref{sssec:moo:nise:w_inf}) of the segment that contains the solutions $\vec{r}^i$ and $\vec{r}^j$ and the intersection point $\underline{\vec{r}}$ between the solution hyperplanes (${\vec{w}^i}^\top \underline{\vec{r}} = {\vec{w}^i}^\top \vec{r}^i$ and ${\vec{w}^j}^\top \underline{\vec{r}} = {\vec{w}^j}^\top \vec{r}^j$) of the neighborhood, thus producing:

\begin{equation}\label{eq:moo:estim_err}
\mu^{i,j}=\sqrt{\frac{(\vec{w}^\top \vec{r}^j-\vec{w}^\top \underline{\vec{r}})^2}{\left|\left|\vec{w}\right|\right|^2}}.
\end{equation}

\begin{figure}[ht!]
\centering
\includegraphics[width=0.55\textwidth]{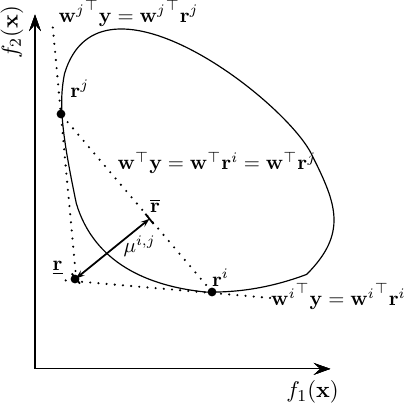}
\caption{Geometrical view of the current inner ($\overline{\vec{r}}$) and outer ($\underline{\vec{r}}$) approximations of the Pareto frontier}
\label{fig:moo:estim_err}
\end{figure}

Figure \ref{fig:moo:estim_err} illustrates the steps involved. Vectors $\vec{w}^i$, $\vec{w}^j$ indicate the weighting vectors used to find the solutions $\vec{r}^i$ and $\vec{r}^j$, respectively. Then, the intersection of ${\vec{w}^{i}}^{\top}\vec{y}={\vec{w}^{i}}^{\top}\vec{r}^i$ and ${\vec{w}^{j}}^{\top}\vec{y}={\vec{w}^{i}}^{\top}\vec{r}^j$, provides $\underline{\vec{r}}$, leading to the distance $\mu^{i,j}$ between $\underline{\vec{r}}$ and $\vec{w}^{\top}\vec{y}=\vec{w}^{\top}\vec{r}^j$ produced by (\ref{eq:moo:estim_err}).

\subsubsection{Calculation of the weighting vector} \label{sssec:moo:nise:w_inf}

Given the neighborhood $(i,j)$ composed of two efficient solutions $\{\vec{r}^i,\vec{r}^j\}$, it is possible to calculate the unitary normal vector $\vec{w}$ of the segment containing these points, by solving the following linear system:

\begin{equation}\label{eq:moo:nise_w_estim}
\begin{cases}
\vec{w}^\top \vec{r}^i = \vec{w}^\top \vec{r}^j\\
\vec{w}^\top \vec{1} = 1.
\end{cases}
\end{equation}

Having determined the weighting vector $\vec{w}$ with maximum error $\mu$ at the $k$-th iteration, it is possible to solve the weighted sum problem (Definition \ref{def:moo:weight_prob}) and find $\vec{r}^k=\vec{f}({\vec{x}^*}^k)$.

\subsubsection{Updating new neighborhoods} \label{sssec:moo:nise:branch}

Given a solution $\vec{r}^k$ associated with the current neighborhood, the new neighborhoods $(i, k)$ and $(k,j)$ are added to $\mathcal{N}^{k}$ and the previous neighbohood $(i, j)$ is deleted, resulting in $\mathcal{N}^{k+1} = \mathcal{N}^{k} \cup \{(i,k),(k,j)\} \setminus (i,j)$.

\subsubsection{Stopping criterion}\label{sssec:moo:nise:stop_crit}
The stopping criterion is satisfied when the largest estimation error $\mu$, defined in Section \ref{sssec:moo:nise:choose}, is smaller than the threshold error $\mu^{stop}$, or a prescribed number of efficient solutions have been generated.

\subsection{Discussion} \label{ssec:moo:nise_disc}

Aiming at further investigating the properties of NISE, we provide graphical examples and theoretical proofs. For convenience, theoretical proofs are presented in Appendices A to D.

\begin{figure}[ht!]
\centering
\subfloat[Initialization]{\includegraphics[width=0.33\textwidth]{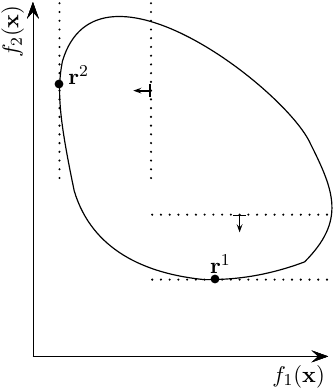}}
\subfloat[First iteration]{\includegraphics[width=0.33\textwidth]{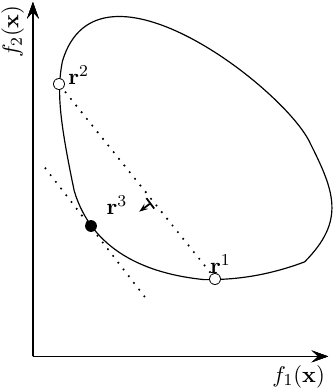}}
\subfloat[Second iteration]{\includegraphics[width=0.33\textwidth]{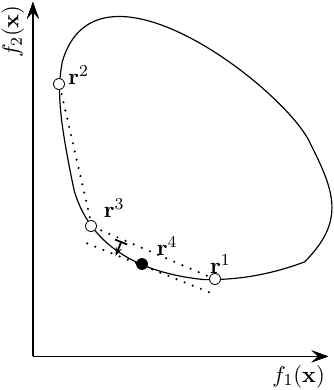}}
 \caption{Illustrative sequence of steps of the NISE method}
 \label{fig:moo:nise}
\end{figure}

A brief illustrative example of the execution of the method is presented in Figure \ref{fig:moo:nise}. The method initiates in Figure \ref{fig:moo:nise}-a, with the determination of the extreme solutions of the problem ($\vec{r}^{1}$ and $\vec{r}^{2}$). In Figure \ref{fig:moo:nise}-b, the unitary normal vector of the segment containing solutions $\vec{r}^1$ and $\vec{r}^2$ is determined, and then Definition \ref{def:moo:weight_prob} is used to find solution $\vec{r}^{3}$. Finally, in Figure \ref{fig:moo:nise}-c, we have two neighborhoods ($\vec{r}^{1}$,$\vec{r}^{3}$) and ($\vec{r}^{3}$,$\vec{r}^{2}$), and the first neighborhood is selected (given the larger margin error $\mu$), thus guiding to solution $\vec{r}^{4}$ using Definition \ref{def:moo:weight_prob} again. This procedure is repeated until convergence, when the maximum margin error considering all the existing neighborhoods is smaller than $\mu^{stop}$.

\begin{figure}[ht!]
\centering
\subfloat[]{\includegraphics[width=0.55\textwidth]{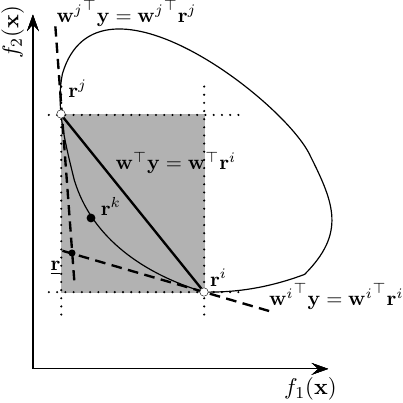}}
\subfloat[]{\includegraphics[width=0.45\textwidth]{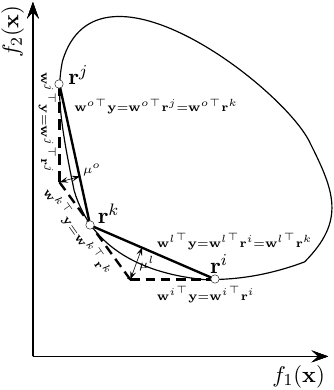}}
 \caption{Geometrical support for the recursivity of the NISE method.}
 \label{fig:moo:nise_an}
\end{figure}

Figure \ref{fig:moo:nise_an}-a represents the main properties of NISE: 1. The weighting vector $\vec{w}$ is located inside the generated cone produced by the weighting vectors $\vec{w}^i$ and $\vec{w}^j$ associated with solutions $\vec{r}^i$ and $\vec{r}^j$ (for a demonstration, see Theorem \ref{theo:moo:weight_2Drecur} of Appendix \ref{appendix:nise_recursivity}); 2. The solution $\vec{r}^k$ found using $\vec{w}$ is inside the gray box delimited by $\vec{r}^i$ and $\vec{r}^j$, as demonstrated in Theorem \ref{theo:moo:weight_2Dlocal} of Appendix \ref{appendix:nise_recursivity}. The first property guarantees an intermediate slope at any step of the algorithm. Given that, both properties guide to a solution between the previous neighborhood solutions. These properties comprise the necessary ingredients to support the recursivity of the method.

Unfortunately, those properties are not retained for problems with more than two objectives: 1. It is not possible to guarantee the nonnegativity of all the elements of the weighting vector $\vec{w}$, thus not satisfying the first property. This proposition is theoretically proved in Theorem \ref{theo:moo:weight_3Dnonrecur} of Appendix \ref{appendix:nise_recursivity}; 2. Even under positivity of the elements of $\vec{w}$, the solution is not necessarily located between the solutions that generate $\vec{w}$. This proposition is theoretically proved in Theorem \ref{theo:moo:weight_3Dnonlocal} of Appendix \ref{appendix:nise_recursivity}.

However, there is another perspective for the approach of the weighted sum method in the context of the NISE algorithm that helped us to construct a more general formulation. In Figure \ref{fig:moo:nise_an}-b we can see the dashed segments associated with the weighting vectors $\vec{w}^i$, $\vec{w}^j$ and $\vec{w}^k$ used to generate the solutions $\vec{r}^i$, $\vec{r}^j$ and $\vec{r}^k$ respectivelly. Due to the optimality of the weighted sum method, it is not possible to find any solution beyond these segments closer to the origin, implying that the dashed segments represent an outer approximation of the frontier. On the other hand, the solid segments in Figure \ref{fig:moo:nise_an}-b, built connecting $\vec{r}^i$ to $\vec{r}^k$ and $\vec{r}^j$ to $\vec{r}^k$, are convex combinations of $\vec{r}^i$ and $\vec{r}^k$ and of $\vec{r}^k$ and $\vec{r}^j$, precluding any solution dominated by these combinations to be found by the weighted sum method, thus creating an inner approximation of the efficient solutions attainable by the weighted sum method.

Given that, we can determine the variables $\underline{\vec{r}}$, $\overline{\vec{r}}$ and $\vec{w}$ which correspond to, respectively, a point limited by the relaxed constraints, a point limited by the inner approximation of the frontier, and a slope for the weighted sum method. Finally, we can define an optimization problem that searches for the maximum distance between the inner approximation $\vec{w}^\top\overline{\vec{r}}$ and the outer approximation $\vec{w}^\top\underline{\vec{r}}$, interpreted here as an extension of the original NISE method and representing the main contribution of this paper.

\section{Extension of the NISE algorithm for more than two objectives} \label{sec:monise}

In this section we present a novel adaptive multi-objective optimization algorithm, an extension of NISE \citep{Cohon1979} capable of dealing with two or more objectives. The main distinct aspect of the proposed methodology is a new optimization model described in Definition \ref{def:moo:monise_w_estim}, responsible for recursively finding the next weighting vector $\vec{w}$ and the current estimation error $\mu$. This extension will be called Many Objective NISE (MONISE) which reduces to a procedure similar to NISE when only two objectives are considered (see Section \ref{sssec:moo:monise:disc} for an explicit comparison).

\subsection{Inner and outer approximation in the context of the weighted sum method}

Consider the utopian solution $\vec{z}^{utopian}$, as well as $L \geq 1$ efficient solutions $\vec{f}(\vec{x}^{i}): i \in \{1,\ldots,L\}$ obtained by the weighted sum method (see Definition \ref{def:moo:weight_prob}) using weighting vectors $\vec{w}^{i}: i \in \{1,\ldots,L\}$. The problem in Definition \ref{def:moo:monise_w_estim} determines a new weighting vector $\vec{w}$, the inner approximation $\overline{\vec{r}}$ and the outer approximation $\underline{\vec{r}}$ of the Pareto frontier associated with the largest distance $\mu$.

The \textbf{outer approximation} is a theoretical limitation for any efficient solution $\vec{x}^*$ attainable by the weighted sum method. So, it is possible to conclude that an outer approximation $\underline{\vec{r}} \in \mathbb{R}^m$ of the objective vector will obey the inequalities ${\vec{w}^{i}}^\top \vec{f}(\vec{x}^{*}) \geq {\vec{w}^{i}}^\top \underline{\vec{r}} \geq {\vec{w}^{i}}^\top \vec{f}(\vec{x}^{i})\ \forall i \in \{1,\ldots,L \}$, since $\vec{x}^{i}$ is the optimal solution of the problem in Definition \ref{def:moo:weight_prob} considering the weighting vector $\vec{w}^{i}$.

The \textbf{inner approximation} is a theoretical limitation for any efficient solution $\vec{x}^*$ attainable by the weighted sum method. Thus there is a weighting vector $\vec{w}$ whose corresponding efficient solution is $\vec{x}^*$, and the inner approximation of the objective vector is $\overline{\vec{r}} \in \mathbb{R}^m$. Following the premises it is possible to conclude that ${\vec{w}}^\top \vec{f}(\vec{x}^{*}) \leq {\vec{w}}^\top \overline{\vec{r}} \leq {\vec{w}}^\top\vec{f}(\vec{x}^i)\ \forall i \in \{1,\ldots,L \}$, since $\vec{x}^*$ is the optimal solution of the problem in Definition \ref{def:moo:weight_prob} considering the weighting vector $\vec{w}$.

Hence, there is an outer estimation $\underline{\vec{r}}$ associated with the outer approximation (${\vec{w}^{i}}^\top \underline{\vec{r}} \geq {\vec{w}^{i}}^\top \vec{f}(\vec{x}^{i})$) and an inner estimation $\overline{\vec{r}}$ associated with the inner approximation (${\vec{w}}^\top\overline{\vec{r}} \leq {\vec{w}}^\top\vec{f}(\vec{x}^i)$), which can be used to define the space of all solutions attainable by the weighted sum method considering the information provided by $L$ previous solutions.

\subsection{Calculating the weighting vector for the weighted sum method} \label{sssec:moo:monise:w_inf}

The calculation of the weighting vector $\vec{w}$ at each iteration of MONISE is done by finding the largest distance between the hyperplanes $\vec{w}^\top \overline{\vec{r}}$ and $\vec{w}^\top \underline{\vec{r}}$.

\begin{definition}\label{def:moo:monise_w_estim}

\begin{equation*}
\begin{aligned}
& \underset{\vec{w}, \overline{\vec{r}}, \underline{\vec{r}}}{\text{maximize}}
& & \mu = \vec{w}^\top \overline{\vec{r}} - \vec{w}^\top \underline{\vec{r}},\\
& \text{subject to}
& & {\vec{w}^{i}}^\top \underline{\vec{r}} \geq {\vec{w}^{i}}^\top \vec{f}(\vec{x}^{i})& \forall i \in \{1,\ldots,L\},\\
& & & {\vec{w}}^\top\overline{\vec{r}} \leq {\vec{w}}^\top\vec{f}(\vec{x}^i)& \forall i \in \{1,\ldots,L\},\\
& & & \underline{\vec{r}} \geq \vec{z}^{utopian},\\
& & & \vec{w} \geq \vec{0},\\
& & & \vec{w}^\top \vec{1} = 1.
\end{aligned}
\end{equation*}
\end{definition}

By solving the problem formulated in Definition \ref{def:moo:monise_w_estim}, the next neighborhood and the weighting vector $\vec{w}$ are determined, as already done by the NISE method. In two dimensions, MONISE searches for an outer approximation $\underline{\vec{r}}$ in the intersection of two slopes as well as finds an inner approximation $\overline{\vec{r}}$ in the convex combination of already found solutions. A mathematical proof can be found in Theorem \ref{theo:moo:monise_nise} of Appendix \ref{appendix:monise_ap}. As we already anticipated, with three or more dimensions, the positivity of each element of the weighting vector $\vec{w}$ is not ensured for the NISE method, which demands novel solution strategies as the one just described for the MONISE method.

Additionally, the problem formalized in Definition \ref{def:moo:monise_w_estim} is bounded. A mathematical proof of this property can be found in Theorem \ref{theo:moo:monise_lim} of Appendix \ref{appendix:monise_ap}.

\subsection{Mixed-integer linear equivalent formulation}

The problem in Definition \ref{def:moo:monise_w_estim} is non-convex, with a bi-convex structure. To solve it we use KKT conditions and some algebraic manipulations (better explained in Appendix \ref{appendix:monise:lip}) to build the equivalent mixed-integer linear formulation in Definition \ref{def:moo:monise_w_estim_pli}. 

The main differences between the formulations in Definitions \ref{def:moo:monise_w_estim} and \ref{def:moo:monise_w_estim_pli} reside in $\vec{w}^\top \overline{\vec{r}}$ being replaced by $v \in \mathbb{R}$ and the addition of the dual variables $\kappa \in \mathbb{R}^L$, $\nu \in \mathbb{R}^m$ and a scalar variable $\mu \in \mathbb{R}$. New equalities derived from KKT conditions are also included, together with mixed-integer linear inequalities representing complementary conditions. These mixed-integer linear inequalities employ binary variables $\kappa^B \in \{0,1\}^L$, $\nu^B \in \{0,1\}^m$.

Using this new formulation, it is possible to solve this problem with well known commercial solvers, such as Gurobi\footnote{Available at \url{www.gurobi.com}} and CPLEX\footnote{Available at \url{www.ibm.com/analytics/cplex-optimizer}}, as well as with open-source softwares such as GLPK\footnote{Available at \url{www.gnu.org/software/glpk}} and CBC\footnote{Available at \url{projects.coin-or.org/Cbc}}.

\begin{definition}\label{def:moo:monise_w_estim_pli}

\begin{equation*}
\begin{aligned}
& \underset{\vec{w}, \overline{\vec{r}}, v, \kappa, \kappa^B, \nu, \nu^B, \mu}{\text{maximize}}
& & \mu\\
& \text{subject to}
& & {\vec{w}^{i}}^\top \underline{\vec{r}} \geq {\vec{w}^{i}}^\top \vec{f}(\vec{x}^{i}),& \forall i \in \{1,\ldots,L\}\\
& & & \underline{\vec{r}} - \sum_{i=1}^L \kappa_i \vec{f}(\vec{x}^{i})  - \nu + \mu \vec{1} = \vec{0},\\
& & & \underline{\vec{r}} \geq \vec{z}^{utopian},\\
& & &  w_i \geq 0,\ w_i \leq \nu_i^B,& \forall i \in \{1,\ldots,m\},\\
& & & \nu_i \geq 0,\ \nu_i \leq (1-\nu_i^B)\overline{\nu_i},& \forall i \in \{1,\ldots,m\},\\
& & & [\vec{w}^\top\vec{f}(\vec{x}^{i})-v] \geq 0,& \forall i \in \{1,\ldots,L\},\\
& & & [\vec{w}^\top\vec{f}(\vec{x}^{i})-v] \leq \kappa_i^B \overline{\kappa_i},& \forall i \in \{1,\ldots,L\},\\
& & & \kappa_i \geq 0,\ \kappa_i \leq (1-\kappa_i^B),& \forall i \in \{1,\ldots,L\},\\
& & & \sum_{i=1}^m w_i = 1,\ \sum_{i=1}^L \kappa_i = 1,
\end{aligned}
\end{equation*}
where $\overline{\kappa_i}$ is a constant upper bound always larger than $[\vec{w}^\top\vec{f}(\vec{x}^{i})-v]$, and $\overline{\nu_i}$ is a constant upper bound always larger than $\nu_i$.
\end{definition}

\subsection{Outline of the methodology}

The Many-Objective NISE method (MONISE), jointly determines the weighting vector and the estimation error. This is done without any additional structure (such as the neighborhood in NISE method), simply resorting to the previous solutions $\{\vec{x}^1, \ldots, \vec{x}^L\}$ and weighting vectors $\{\vec{w}^1, \ldots, \vec{w}^L\}$. The procedure adopted by MONISE is simpler and more general than the one required by NISE and may be summarized in three phases: (1) The \textbf{initialization} (described in Section \ref{sssec:moo:monise:init}); (2) The \textbf{iterative  process} responsible for defining, at the $L$-th iteration, the weighting vector $\vec{w}^{L+1}$ (described in Section \ref{sssec:moo:monise:choose}) used to find the solution $\vec{x}^{L+1}$ when optimizing the problem in Definition \ref{def:moo:weight_prob}; and finally, (3) The \textbf{stopping criterion} described in Section \ref{sssec:moo:monise:stop_crit}.

\subsubsection{Initialization} \label{sssec:moo:monise:init}

The algorithm is initialized with the following inputs: (1) the utopian solution $\vec{z}^{utopian}$, and (2) at least one weighting vector ($\vec{w}^i$) and its corresponding solution ($\vec{x}^i$). 

This input information is enough to guarantee that problem in Definition \ref{def:moo:monise_w_estim} is limited. However, since $\vec{z}^{utopian}$ is obtained using the $k$-th individual minimum values, we included all $k$-th individual minimum solutions and their weighting vectors (null values for all objectives but the $k$-th objective) in the initialization, thus producing $m$ minimum solutions.

\subsubsection{Choice of the next weighting vector for the weighted sum method} \label{sssec:moo:monise:choose}

The iterative choice of the weighting vector is given by solving problem in Definition \ref{def:moo:monise_w_estim}, recursively using the outcomes of the weighted sum problems already solved. The approximation error ($\mu$) of the iteration is also determined. In practice, the equivalent problem in Definition \ref{def:moo:monise_w_estim_pli} is solved instead of the one in Definition \ref{def:moo:monise_w_estim}.

\subsubsection{Stopping criteria} \label{sssec:moo:monise:stop_crit}

The execution of MONISE stops when the estimation error $\mu$ is lower than a threshold $\mu^{stop}$ or when the number of efficient solutions generated by the method reaches a pre-specified value.

\subsection{Discussion} \label{sssec:moo:monise:disc}

\begin{figure}[ht!]
\centering
\subfloat[]{\includegraphics[width=0.35\textwidth]{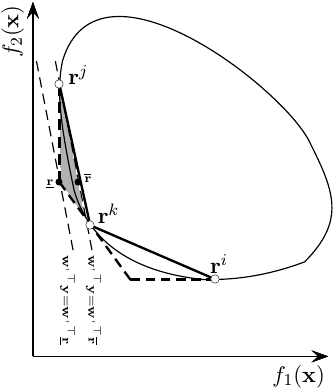}}
\subfloat[]{\includegraphics[width=0.35\textwidth]{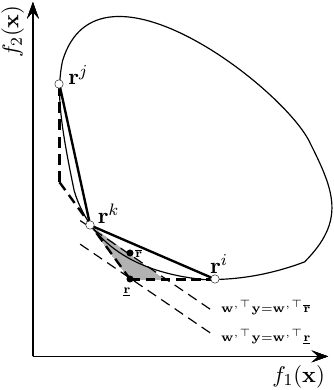}}\\
\subfloat[]{\includegraphics[width=0.35\textwidth]{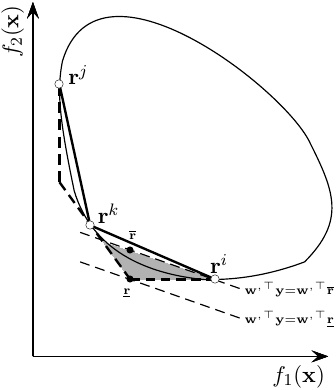}}
\subfloat[]{\includegraphics[width=0.35\textwidth]{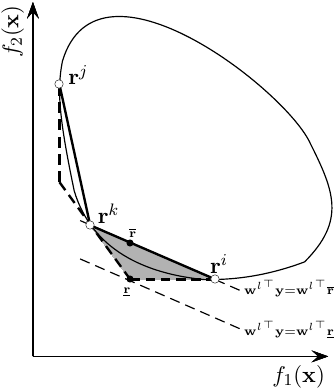}}
 \caption{Sequence of suboptimal steps until achieving the optimal solution of the weighting vector calculation (described in Definitions \ref{def:moo:monise_w_estim} and \ref{def:moo:monise_w_estim_pli}) of the MONISE method.}
 \label{fig:moo:monise_an}
\end{figure}

\begin{figure}[ht!]
\centering
\subfloat[First iteration]{\includegraphics[width=0.33\textwidth]{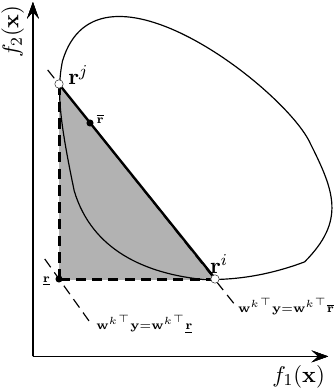}}
\subfloat[Second iteration]{\includegraphics[width=0.33\textwidth]{fig/monise_2.pdf}}
\subfloat[Third iteration]{\includegraphics[width=0.33\textwidth]{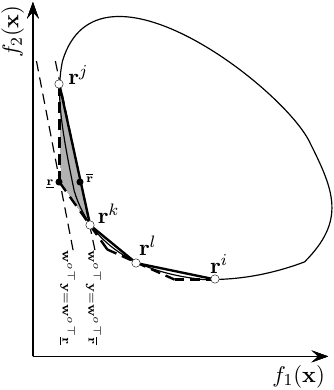}}
 \caption{Illustrative sequence of optimal steps of the MONISE method}
 \label{fig:moo:monise}
\end{figure}

The main properties of MONISE will be illustrated by some graphical examples. Figures \ref{fig:moo:monise_an} and \ref{fig:moo:monise} are restricted to a two-objectives scenario for visualization purposes.
The optimization steps based on the problems in Definitions \ref{def:moo:monise_w_estim} and \ref{def:moo:monise_w_estim_pli} are represented in Figure \ref{fig:moo:monise_an}. Figures \ref{fig:moo:monise_an}-a to \ref{fig:moo:monise_an}-c show sub-optimal weighting vectors $\vec{w}^{,}$, the black dots are the optimized $\overline{\vec{r}}$ and $\underline{\vec{r}}$ when $\vec{w}^{,}$ is considered, and the gray area is the feasible region for these variables. It is possible to see the margin ${\vec{w}^{,}}^\top\overline{\vec{r}}-{\vec{w}^{,}}^\top\underline{\vec{r}}$ gradually increasing from the first to the third illustration. It begins in a suboptimal ``neighborhood'' and increases until it finds $\vec{w}^{,}$. That allows $\overline{\vec{r}}$ to become closer to the segment that connects $\vec{r}^k$ and $\vec{r}^i$, creating a continuous behavior which ends up with the weight vector chosen by the NISE method. However, instead of working with the concept of neighborhood, all the regions not constrained by the optimality of the weighted sum method are explored. The mixed-integer linear optimization process enables this automatic search with global properties, allowing the process to find the optimal value observed in Figure \ref{fig:moo:monise_an}-d. Essentially, this is what gives MONISE the ability to deal with two or more objectives.

The optimization procedure of Definition \ref{def:moo:monise_w_estim}, and its equivalent formulation of Definition \ref{def:moo:monise_w_estim_pli}, can be seen as a search for $\vec{w}$ that ultimately guides to the maximum margin ${\vec{w}}^\top\overline{\vec{r}}-{\vec{w}}^\top\underline{\vec{r}}$ considering $\overline{\vec{r}}$ and $\underline{\vec{r}}$ constrained by the optimality of the weighted sum method. Figure \ref{fig:moo:monise} depicts this procedure in distinct scenarios, starting with the first iteration after initialization, and showing the optimal $\vec{w}$ at each iteration. 
Given that we have only two objectives here, MONISE is capable of automatically finding a valid neighborhood for the next iteration using an iterative procedure similar to NISE. The only distinction comes from the margin calculation which is $\sqrt{\frac{(\vec{w}^\top \vec{r}^j-\vec{w}^\top \underline{\vec{r}})^2}{\left|\left|\vec{w}\right|\right|^2}}$ (where $\vec{r}^j\equiv \overline{\vec{r}}$) for NISE and ${\vec{w}}^\top\overline{\vec{r}}-{\vec{w}}^\top\underline{\vec{r}}$ for MONISE. This subtle difference may sometimes result in distinct neighborhood choices. We could adopt the MONISE's margin calculation to make NISE fully equivalent to MONISE. On the other hand, we could not use the NISE's margin calculation in the objective function of the MONISE formulation in Definition \ref{def:moo:monise_w_estim}, because it would preclude the existence of a formulation equivalent to the one in Definition \ref{def:moo:monise_w_estim_pli}. Notice that this reasoning only makes sense in the case of two objectives.

Now, let us consider a case study with three objectives.

\begin{definition}\label{eq:moo:cex-13D}
\begin{equation*} 
\begin{aligned}
& \underset{\vec{x}}{\text{minimize}}
& & \vec{f}(\vec{x}) = \left[(x_1-1)^2,(x_2-1)^2,(x_3-1)^2\right]\\
& \text{subject to}
& & \vec{x}^\top \vec{1} = 1,\ x_1, x_2, x_3 \geq 0.
\end{aligned}
\end{equation*}
\end{definition}

\begin{figure}[ht!]
\centering
\subfloat{\includegraphics[width=0.5\textwidth]{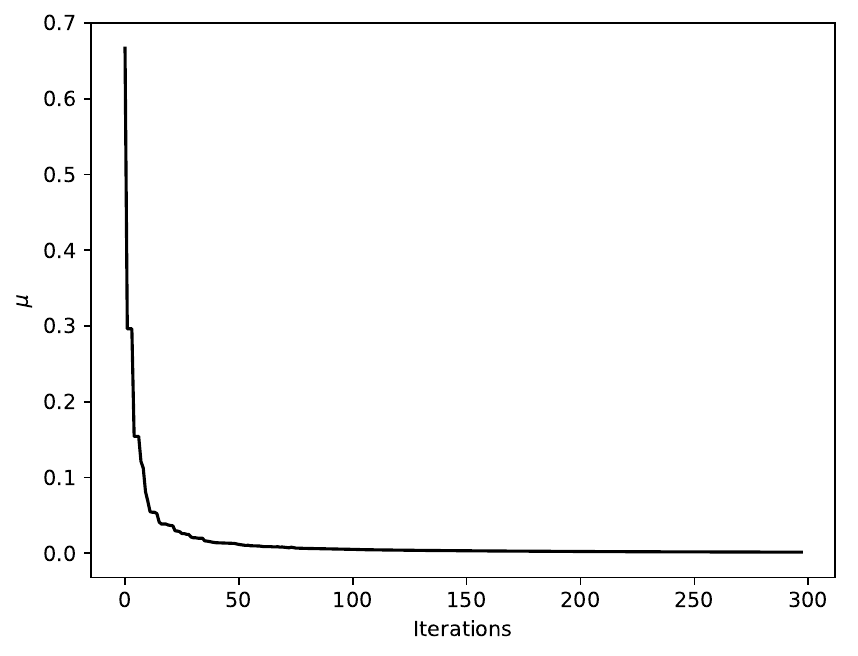}}
\subfloat{\includegraphics[width=0.5\textwidth]{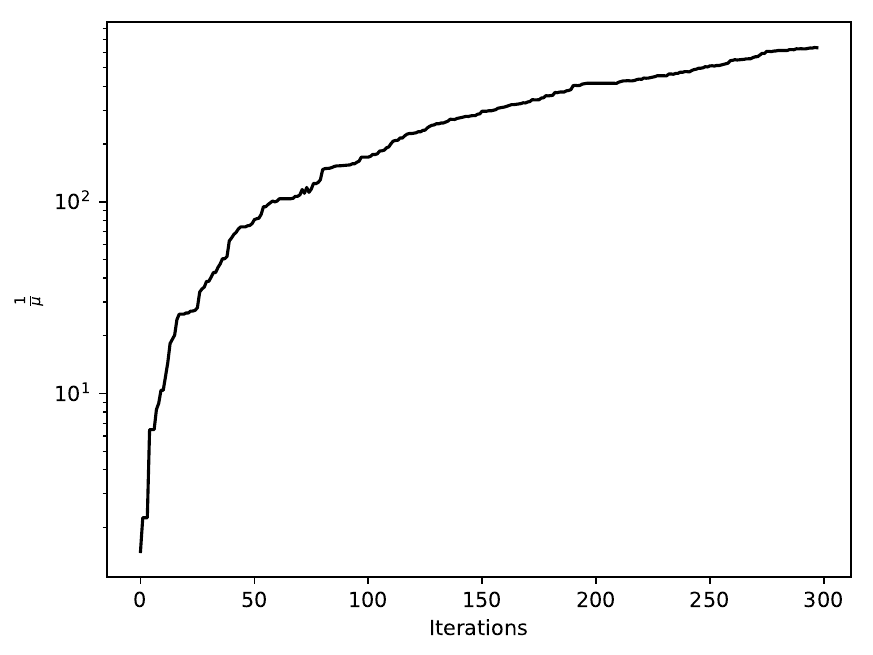}}
 \caption{Evolution in linear (for $\mu$) and logarithmic (for $\frac{1}{\mu}$) scale of the margin $\mu$ along iterations for the problem in Definition \ref{eq:moo:cex-13D}.}
 \label{fig:moo:monise_margin}
\end{figure}

\begin{figure}[ht!]
\centering
\subfloat{\includegraphics[width=0.5\textwidth]{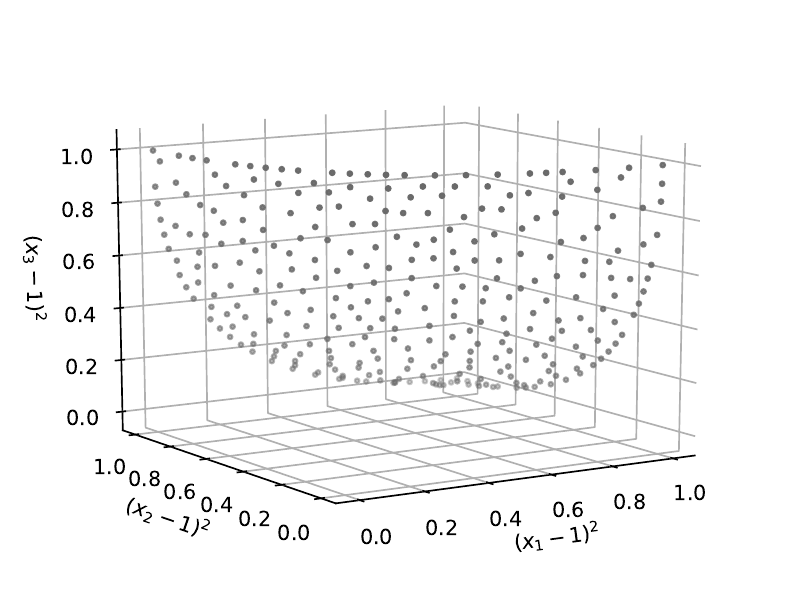}}
\subfloat{\includegraphics[width=0.5\textwidth]{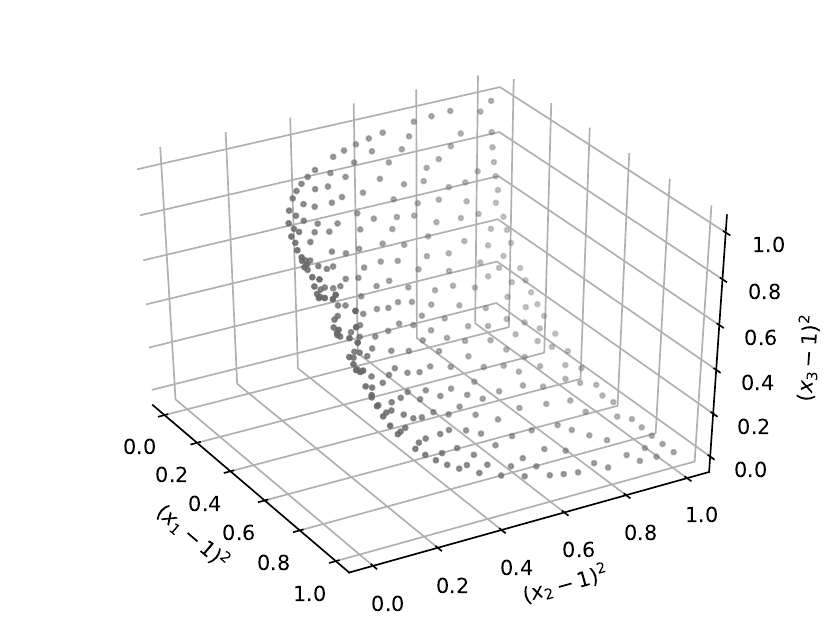}}
 \caption{Two perspectives of the non-inferior set automatically obtained at the Pareto frontier for the problem in Definition \ref{eq:moo:cex-13D} using MONISE.}
 \label{fig:moo:monise_front}
\end{figure}

The simple problem presented in Definition \ref{eq:moo:cex-13D} will be used to further investigate the behavior of MONISE. In Figure \ref{fig:moo:monise_margin} it is shown the evolution of the optimized margin $\mu$ along the iterations, which monotonically decreases to zero in few iterations. Figure \ref{fig:moo:monise_front} shows a well distributed sample of the Pareto frontier after 300 iterations in two perspectives.

A more realistic and comparative experimental analysis is presented in Section \ref{sec:experiments}.

\section{Experiments} \label{sec:experiments}
  
The proposed methodology is a multi-objective method that iteratively finds efficient solutions using the weighted sum method. Its adaptive procedure uses mixed-integer linear programming to automatically define, in a constrained objective space, the weighting vector that best improves the Pareto frontier representation. This constrained objective space is built excluding any candidate solution that contradicts optimality premises of the solutions found by the weighted sum method. Since this method finds each efficient solution using a weighted sum method, the objective function at each iteration consists in a linear combination of the objectives, not increasing the level of complexity of the optimization problem (for instance, there is no need of any additional structural information such as the definition of a neighborhood). 

These features characterize the proposed method as a simple, interpretable and appealing multi-objective method, that, at each iteration, globally and automatically explores the Pareto frontier to recursively find new efficient solutions. Being based on scalarization, MONISE is flexible enough to admit any global solver to the single-objective problem that is produced every time that the next weighting vector is automatically determined after solving the mixed-integer linear problem in Definition 9. Consequently, the solvability of this single-objective problem obtained after determining the weighting vector is restricted only by the existence of a corresponding global solver.

To validate this algorithm, we consider two problems as benchmarks: (1) a multilabel classification, which is a nonlinear convex problem (solved using optimization tools in the SciPy Python library\footnote{Available in \url{https://www.scipy.org/}}), and (2) the knapsack problem, which is a combinatorial problem (solved using optimization tools in the Gurobi Python library\footnote{Available at \url{http://www.gurobi.com/}}). The metrics used to evaluate the performance of the methods are: (1) hypervolume, and (2) execution time. Multilabel classification problems are generally characterized by convex Pareto frontiers, which are properly explored by MONISE. However, every multi-objective approach founded on scalarization, such as NISE and MONISE, will not be able to guarantee the exploration of the whole Pareto frontier when convexity is no more assured, because efficient solutions that are dominated by convex combinations of other efficient solutions are not achievable by scalarization. Therefore the application of MONISE to deal with discrete Pareto frontiers, such as the ones characterizing knapsack problems, is challenging and motivates the experiments.

The definition of the multilabel classification and the knapsack problems will be presented in Sections \ref{ssec:exp:mlc} and \ref{ssec:exp:ksp}, respectively, and the hypervolume metric is briefly reviewed in Section \ref{ssec:exp:eval}.

\subsection{Multilabel classification} \label{ssec:exp:mlc}

Multilabel classification is a relevant problem on supervised learning: given a sample $\vec{x}_{i}$, the objective is to find a subset of labels related to this sample, thus admitting more than one label per sample.

It is supposed that there exist $m$ labels and $N$ samples: $\vec{x}_{i} \in \mathbb{R}^d,\ i \in \{1,\ldots,N\}$, represents the input feature vector and $y_{i}^{(l)} \in \{0,1\},\ l \in \{1,\ldots,m\},\ i \in \{1, \ldots, N\}$, is the membership of sample $i$ to label $l$, which is the value that we want to predict. The decision variables are the components of the vector $\theta  \in \mathbb{R}^{n} \equiv \mathbb{R}^{d+1}$. It is then possible to devise the following multi-objective model:

\begin{equation} \label{eq:ml:mtl}
\begin{aligned}
& \underset{\theta}{\text{minimize}}
& & \vec{f}(\vec{x})\equiv \{f_{(1)}(\theta, \vec{x}, \vec{y}), \ldots, f_{(m)}(\theta, \vec{x}, \vec{y})\}
\end{aligned}
\end{equation}
where $f_{(l)}(\theta) =  \sum_{i=1}^{n} -\left[ y_i^{(l)} \ln \left(\frac{e^{\theta^\top \phi(\vec{x}_i)}}{1+e^{\theta^\top \phi(\vec{x}_i)}}\right) + (1-y_i^{(l)}) \ln \left(1-\frac{e^{\theta^\top \phi(\vec{x}_i)}}{1+e^{\theta^\top \phi(\vec{x}_i)}}\right) \right]$, $\forall\ l \in \{1,\ldots,m\}$, is the classification loss for label $l$, and $\phi_0(\vec{x}_i)=1, \phi_j(\vec{x}_i)=x_{i,j}\ \forall\ j \not = 0,\ j \in \{1,\ldots,d\}$.

To conduct the test for multilabel classification, we consider five datasets\footnote{Available at \url{mulan.sourceforge.net/datasets-mlc.html}} whose main characteristics are described in Table \ref{tab:dataset}.

\begin{table}[ht!]
\centering
\caption{Description of the main characteristics of the five multilabel datasets.}
\label{tab:dataset}
\begin{tabular}{l|r|r|r}
name     & instances $(N)$ & features $(d)$ & labels $(m)$ \\ \hline
emotions & 593  & 72   & 6   \\
flags    & 194  & 19   & 7   \\
yeast    & 2,417 & 103  & 14  \\
birds    & 645  & 260  & 19  \\
genbase  & 662  & 1,186 & 27 
\end{tabular}
\end{table}

\subsection{The knapsack problem} \label{ssec:exp:ksp}

This problem of central practical and theoretical relevance in linear-integer programming, leads to discrete and non-convex Pareto frontiers \citep{Bazgan2009a}, and imposes a challenge for multi-objective optimization methods. As an additional motivation, methods for constructing controled knapsack instances are available in the literature \citep{Bazgan2009a}.

Given a constraint of capacity $T$, one must choose between $n$ items with diversity of sizes and utility values to fill up the knapsack. The goal is to maximize the utility value of the picked items without exceeding the capacity of the knapsack. In the multi-objective case, each item $i$ with size $t_i$ has $m$ distinct utility values $v_i^1,v_i^2,\ldots,v_i^m$ and the idea is to maximize all these utility values concurrently.

The vector of decision variables $\vec{x}$ of the problem is taken as a binary vector. The interpretation of the vector indicates if an item was picked ($x_i=1$) or not ($x_i=0$).

\begin{equation}
\begin{aligned}
& \underset{\vec{x}}{\text{minimize}}
& & \vec{f}(\vec{x})\equiv \{f_{(1)}(\vec{x}), \ldots, f_{(m)}(\vec{x})\}\\
& \text{where}
& & f_{(l)}(\vec{x}) = \sum_{i=1}^{n}v_{i}^{l}x_{i}, \forall l \in \{1,2,\ldots,m \},\\
& \text{subject to}
& & \sum_{i=1}^{n}t_{i}x_{i}\leq T,\\
& & & x_{i}\in\{0,1\}, \forall i \in \{1,\ldots,n\}.\\
\end{aligned}
\end{equation}

The generation of the instances was made using a procedure suggested by \citet{Bazgan2009a}. All parameters of the knapsack problem were randomly generated inside a given interval, forcing a conflicting behavior between the multiple objectives. Each size $t_i$ and values $v_i^1,v_i^2,\ldots,v_i^m$ are generated in the interval $[0,1000]$. In this experiment we use 100 items and capacities of $T=1000$ and $T=2000$. The tighter constraint ($T=1000$) induces a harder conflict among the objectives, so that the few items in the obtained efficient solutions have a high chance of displaying distinct values, while the looser constraint ($T=2000$) leads to scattered solutions. Any upper capacity bound greater than $T=2000$ would lead to solutions with too many items in the knapsack, reducing the conflict between the objectives and its practical relevance.

\subsection{Hypervolume - Evaluation metric} \label{ssec:exp:eval}

The hypervolume metric \citep{Fleischer2003} may be described as follows: given a reference point that is dominated by all efficient solutions, preferably the Nadir point, it is calculated the hypervolume formed using this point and all previous solutions as delimiters. The hypervolume metric may be obtained in a cumulative manner. The additional volume consists in: considering the set of solutions $R=\{\vec{r}^1,\vec{r}^2,\ldots,\vec{r}^{k-1} \}$ with the hypervolume already calculated $\mathcal{V}^{current}$, the additional volume implied by a solution $\vec{r}^k$ is given by: $\mathcal{V}(\vec{r}^k|R)=\mathcal{V}(\vec{r}^k)-(\mathcal{V}^{current}\cap \mathcal{V}(\vec{r}^k))$.

\begin{figure}[!ht]
\centering
\includegraphics[scale=0.95]{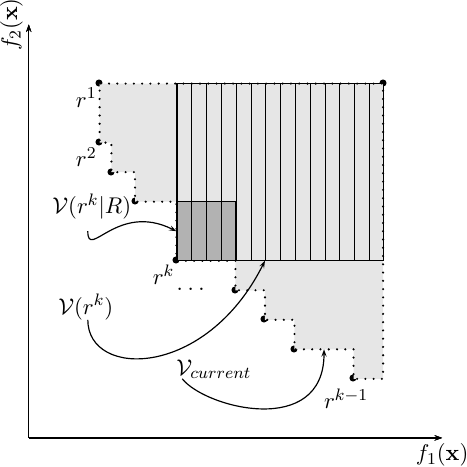}
 \caption{Representation of the hypervolume metric and its recursive updating, when only two objectives are considered}
 \label{fig:moo:hypervolume}
\end{figure}

This metric exhibits interesting properties because dominated solutions do not contribute to the hypervolume. Furthermore, when two solutions are too close, the contribution of the second when the first was already considered would add too little to the hypervolume. In Figure \ref{fig:moo:hypervolume}, the grey volume indicates $\mathcal{V}^{current}$, the shaded volume corresponds to $\mathcal{V}(\vec{r}^k)$ and the dark grey volume represents $\mathcal{V}(\vec{r}^k|R)$.

\subsection{Experimental setup}

The proposed algorithm was compared with four state-of-the-art evolutionary approaches, NSGA-II \citep{Deb2002}, NSGA-III \citep{Deb2013}, SPEA2 \citep{Zitzler2001}, and SMS-EMOA \citep{Beume2007}\footnote{First three algorithms are available at \url{https://github.com/Project-Platypus/Platypus} and the last one is available at \url{https://github.com/marcosmrai/moo_experiments}}. It was also compared with one random approach that uses the Kraemer Algorithm \citep{Smith2004} to randomly sample weight vectors of the weighted sum method -- called here as RAND, and three non-heuristic algorithms, the PGEN \citep{Craft2006}, the proposal of \cite{Rennen2009} -- called here as RENN, and the NC \citep{Messac2003}\footnote{The random, and the three non-heuristic methods are available at \url{https://github.com/marcosmrai/moopt}}. To create the reference points for NSGA-III, we define three uniform partitions of every objective. Since the evolutionary algorithms use a fixed number of solutions, the experiment was designed to look for $5\times m$ efficient solutions, where $m$ is the number of objectives. Since the NC algorithm does not fix the number of solutions, its parameters were fine-tuned until we find a number of solutions as close as possible but not smaller than $5\times m$ (the reported hypervolume comes from the $5\times m$ solutions with the greater contribution to the hypervolume). These eight algorithms were chosen because they are good representatives of methods which are not constrained to any type of optimization problem and accept a large number of objectives. Every iteration of any algorithm is allowed to execute if the run time is up to 3600 seconds, limiting the number of generations in the evolutionary approaches and reducing the number of solutions in the non-heuristic approaches.

For the multilabel classification, the evolutionary algorithms had the following configuration: population of $5\times m$ where $m$ is the number of objectives, with the same amount of individuals being generated in the offspring (except for SMS-EMOA that generates a single individual), 90\% of offspring generated by uniform crossover, 10\% of offspring generated by Gaussian mutation, and the number of generations is enough to runtime be at least 3600 seconds.
For the knapsack problem: population of $5\times m$, with the same amount of individual being generated in the offspring (except for SMS-EMOA that generates a single individual), 90\% of offspring generated by uniform crossover, 10\% of offspring generated by a mutation that can include or exclude an item from the knapsack, and the number of generations is enough to runtime be at least 3600 seconds. The evaluation of the knapsack problem considered the correction strategy \citep{Michalewicz1994}, removing items with the lowest profit ratio (where the profit ratio of an item $i$ is defined by $\frac{\prod_{k=1}^m v_i^k}{t_i}$). In addition, the population had four times the desired number of solutions to increase the evolution pressure, since preliminary experiments showed a lack of performance improvement with a smaller population size.

Since MONISE is based on a mixed-integer linear formulation to calculate the weight vectors, we explored the options of commercial solvers to control the precision and time of the calculations of this step. We used a limit of 2 seconds for each calculation, as well as a gap of 0.01. Worth mentioning that the entire experimental setup, including datasets, is available at \url{github.com/marcosmrai/moo}.

\subsection{Obtained Results}

The experiments were evaluated in terms of execution time and final hypervolume. Aiming at a fair evaluation of the hypervolume, the reference point was chosen by selecting the worst value of each objective from all non-dominated solutions and considering all the nine evaluated algorithms.

Table \ref{table:convex_hv} presents the hypervolume for the multilabel classification instances, Table \ref{table:convex_time} the execution time for the same instances, and Table \ref{table:convex_iter} the number of iterations performed by each algorithm. Table \ref{table:ks_hv} presents the hypervolume for the knapsack problem instances, Table \ref{table:ks_time} the execution time for the same instances, and Table \ref{table:ks_iter} the number of iterations performed by each algorithm. The first columns of Tables \ref{table:convex_hv}, \ref{table:convex_time} and \ref{table:convex_iter} indicate the name of the instances, and the second columns, the number of objectives. The first columns of Tables \ref{table:ks_hv}, \ref{table:ks_time} and \ref{table:ks_iter} indicate the capacity of the knapsack for the instances, and the second columns, the number of objectives. The instances in which there was segmentation fault or memory overload for the method, RENN receives the symbol ``--''. Those errors were caused by the necessity of manipulating a large number of points associated with a high dimensional convex hull.

\ifx
\begin{table}[ht!]
\setlength{\tabcolsep}{3pt}
\centering
\caption{Comparison in terms of hypervolume for the multilabel classification problem.}
\label{table:convex_hv}
\begin{tabular}{r!{\vrule width 1.5pt}r!{\vrule width 1.5pt}r|r|r|r|r|r|r|r|r}
 & $m$ & MONISE   & RAND   & NC       & PGEN     & RENN   & NSGA-II   & NSGA-III  & SPEA2    & SMS-E       \\
\noalign{\hrule height 1.5pt}
emotions & 6  & 0.725 & 0.688 & 0.565 & 0.644 & 0.729 & 0.586 & 0.463 & 0.565 & 0.480 \\
flags    & 7  & 0.796 & 0.754 & 0.773 & 0.711 & 0.802 & 0.724 & 0.551 & 0.708 & 0.544 \\
yeast    & 14 & 0.556 & 0.526 & 0.509 & 0.453 &  --   & 0.383 & 0.322 & 0.136 & 0.101 \\
birds    & 19 & 0.980 & 0.966 & 0.854 & 0.958 &  --   & 0.841 & 0.947 & 0.903 & 0.904 \\
genbase  & 27 & 0.786 & 0.740 & 0.363 & 0.631 &  --   & 0.517 & 0.610 & 0.500 & 0.526
\end{tabular}
\end{table}
\fi 

\begin{table}[ht!]
\setlength{\tabcolsep}{3pt}
\centering
\caption{Comparison in terms of hypervolume for the multilabel classification problem.}
\label{table:convex_hv}
\begin{tabular}{r!{\vrule width 1.5pt}r!{\vrule width 1.5pt}r|r|r|r|r|r|r|r|r}
 & $m$ & MONISE   & RAND   & NC       & PGEN     & RENN   & NSGA-II   & NSGA-III  & SPEA2    & SMS-E       \\
\noalign{\hrule height 1.5pt}
emotions & 6  & \cgs0.725 & \cgt0.696 & 0.565 & 0.644 & \cgf0.730  & 0.589 & 0.459 & 0.544 & 0.457 \\
flags    & 7  & \cgs0.796 & \cgt0.776 & 0.773 & 0.711 & \cgf0.797  & 0.736 & 0.538 & 0.702 & 0.538 \\
yeast    & 14 & \cgf0.532 & \cgs0.523 & \cgt0.509 & 0.447 & --     & 0.135 & 0.398 & 0.112 & 0.256 \\
birds    & 19 & \cgf0.974 & \cgs0.966 & 0.852 & \cgt0.958 & --     & 0.840 & 0.943 & 0.845 & 0.833 \\
genbase  & 27 & \cgf0.776 & \cgs0.741 & 0.359 & 0.631 & --     & \cgt0.674 & 0.108 & 0.437 & 0.461
\end{tabular}
\end{table}

\ifx
\begin{table}[ht!]
\setlength{\tabcolsep}{3pt}
\centering
\caption{Comparison in terms of execution time (in seconds) for the multilabel classification problem.}
\label{table:convex_time}
\begin{tabular}{r!{\vrule width 1.5pt}r!{\vrule width 1.5pt}r|r|r|r|r|r|r|r|r}
$T$ & $m$ & MONISE   & RAND   & NC       & PGEN     & RENN   & NSGA-II   & NSGA-III  & SPEA2    & SMS-E       \\
\noalign{\hrule height 1.5pt}
emotions & 6  & 22  & 2.1 & 39   & 6    & 237  & 3600 & 3600 & 3600 & 3600 \\
flags    & 7  & 33  & 0.4 & 11   & 34   & 1092 & 3600 & 3601 & 3600 & 3600 \\
yeast    & 14 & 113 & 5.6 & 344  & 4461 &  --  & 3601 & 3613 & 3603 & 3603 \\
birds    & 19 & 174 & 4.9 & 3641 & 4438 &  --  & 3603 & 3671 & 3603 & 3601 \\
genbase  & 27 & 249 & 9.9 & 3609 & 4978 &  --  & 3602 & 4097 & 3606 & 3821
\end{tabular}
\end{table}
\fi

\begin{table}[ht!]
\setlength{\tabcolsep}{3pt}
\centering
\caption{Comparison in terms of execution time (in seconds) for the multilabel classification problem.}
\label{table:convex_time}
\begin{tabular}{r!{\vrule width 1.5pt}r!{\vrule width 1.5pt}r|r|r|r|r|r|r|r|r}
$T$ & $m$ & MONISE   & RAND   & NC       & PGEN     & RENN   & NSGA-II   & NSGA-III  & SPEA2    & SMS-E       \\
\noalign{\hrule height 1.5pt}
emotions & 6  & \cgt11  & \cgf2  & 51   & \cgs8    & 716  & 3601 & 3600 & 3601 & 3600 \\
flags    & 7  & \cgt17  & \cgf1  & \cgs14   & 43   & 3968 & 3600 & 3601 & 3601 & 3600 \\
yeast    & 14 & \cgs85  & \cgf9  & \cgt318  & 3805 & --   & 3602 & 3679 & 3605 & 3601 \\
birds    & 19 & \cgs123 & \cgf16 & 3978 & 4278 & --   & 3602 & 3644 & 3603 & 3605 \\
genbase  & 27 & \cgs292 & \cgf42 & 3615 & 4520 & --   & 3607 & 4023 & 3605 & 3820 \\
\end{tabular}
\end{table}

\begin{table}[ht!]
\setlength{\tabcolsep}{3pt}
\centering
\caption{Comparison in terms of number of iterations (number of solutions for the deterministic algorithms and generations for the evolutionary) for the multilabel classification problem.}
\label{table:convex_iter}
\begin{tabular}{r!{\vrule width 1.5pt}r!{\vrule width 1.5pt}r|r|r|r|r|r|r|r|r}
$T$ & $m$ & MONISE   & RAND   & NC       & PGEN     & RENN   & NSGA-II   & NSGA-III  & SPEA2    & SMS-E       \\
\noalign{\hrule height 1.5pt}
emotions & 6  & 30  & 30  & 81  & 30 & 30 & 24264 & 8413 & 13548 & 275042 \\
flags    & 7  & 35  & 35  & 92  & 35 & 32 & 49891 & 5003 & 13823 & 385850 \\
yeast    & 14 & 70  & 70  & 103 & 30 & 0  & 4220  & 93   & 1555  & 5678   \\
birds    & 19 & 95  & 95  & 39  & 28 & 0  & 2148  & 25   & 2283  & 6847   \\
genbase  & 27 & 135 & 135 & 101 & 33 & 0  & 1145  & 5    & 831   & 1336   
\end{tabular}
\end{table}

\ifx
\begin{table}[ht!]
\setlength{\tabcolsep}{5pt}
\centering
\caption{Comparison in terms of hypervolume for the multilabel classification problem.}
\label{table:convex_hv}
\begin{tabular}{l|r|rrrrr}
\hline
                &  $m$ & MONISE & NC    & NSGA-II & NSGA-III & SPEA2 \\ \hline
emotions & 7  & 0.537 & 0.466 & 0.523 & 0.364 & 0.501 \\
flags    & 8  & 0.593 & 0.652 & 0.667 & 0.527 & 0.625 \\
yeast    & 15 & 0.705 & 0.665 & 0.334 & 0.181 & 0.410 \\
birds    & 20 & 0.996 & 0.989 & 0.983 & 0.986 & 0.987 \\
genbase  & 28 & 0.721 & 0.755 & 0.205 & 0.041 & 0.118 \\ \hline
\end{tabular}
\end{table}

\begin{table}[ht!]
\setlength{\tabcolsep}{5pt}
\centering
\caption{Comparison in terms of execution time (in seconds) for the multilabel classification problem.}
\label{table:convex_time}
\begin{tabular}{l|r|rrrrr}
\hline
                   &  $m$ & MONISE & NC    & NSGA-II & NSGA-III & SPEA2 \\ \hline
emotions & 7  & 23   & 37   & 118  & 207   & 189  \\
flags    & 8  & 11   & 61   & 43   & 234   & 99   \\
yeast    & 15 & 556  & 2508 & 1955 & 3418  & 2093 \\
birds    & 20 & 527  & 1147 & 1247 & 3897  & 2703 \\
genbase  & 28 & 1722 & 9477 & 1890 & 15414 & 3435  \\ \hline
\end{tabular}
\end{table}
\fi

\ifx
\begin{table}[ht!]
\setlength{\tabcolsep}{3pt}
\centering
\caption{Comparison in terms of hypervolume for the knapsack problem.}
\label{table:ks_hv}
\begin{tabular}{r!{\vrule width 1.5pt}r!{\vrule width 1.5pt}r|r|r|r|r|r|r|r|r}
$T$ & $m$ & MONISE   & RAND   & NC       & PGEN     & RENN   & NSGA-II   & NSGA-III  & SPEA2    & SMS-E       \\
\noalign{\hrule height 1.5pt}
$1000$ & $5 $ & $5.24\cdot 10^{\text{-}1}$ & $4.8\cdot 10^{\text{-}1}$ & $3.5\cdot 10^{\text{-}1}$ & $4.4\cdot 10^{\text{-}1}$ & $5.26\cdot 10^{\text{-}1}$ & $4.3\cdot 10^{\text{-}1}$ & $4.1\cdot 10^{\text{-}1}$ & $4.5\cdot 10^{\text{-}5}$ & $1.5\cdot 10^{\text{-}2}$ \\
$2000$ & $5 $ & $4.17\cdot 10^{\text{-}1}$ & $4.1\cdot 10^{\text{-}1}$ & $4.0\cdot 10^{\text{-}1}$ & $3.9\cdot 10^{\text{-}1}$ & $4.16\cdot 10^{\text{-}1}$ & $3.2\cdot 10^{\text{-}1}$ & $3.7\cdot 10^{\text{-}1}$ & $4.9\cdot 10^{\text{-}2}$ & $3.3\cdot 10^{\text{-}5}$ \\
$1000$ & $10$ & $1.37\cdot 10^{\text{-}1}$ & $1.2\cdot 10^{\text{-}1}$ & $7.0\cdot 10^{\text{-}2}$ & $1.1\cdot 10^{\text{-}1}$ & $9.71\cdot 10^{\text{-}2}$ & $8.7\cdot 10^{\text{-}2}$ & $1.7\cdot 10^{\text{-}2}$ & $4.9\cdot 10^{\text{-}5}$ & $1.1\cdot 10^{\text{-}3}$ \\
$2000$ & $10$ & $8.42\cdot 10^{\text{-}2}$ & $7.9\cdot 10^{\text{-}2}$ & $2.8\cdot 10^{\text{-}2}$ & $7.0\cdot 10^{\text{-}2}$ & $5.24\cdot 10^{\text{-}2}$ & $5.3\cdot 10^{\text{-}2}$ & $5.7\cdot 10^{\text{-}3}$ & $1.3\cdot 10^{\text{-}3}$ & $1.5\cdot 10^{\text{-}3}$ \\
$1000$ & $15$ & $2.74\cdot 10^{\text{-}2}$ & $2.2\cdot 10^{\text{-}2}$ & $6.5\cdot 10^{\text{-}3}$ & $1.7\cdot 10^{\text{-}2}$ & -- & $1.4\cdot 10^{\text{-}2}$ & $4.0\cdot 10^{\text{-}4}$ & $3.9\cdot 10^{\text{-}7}$ & $4.0\cdot 10^{\text{-}5}$ \\
$2000$ & $15$ & $1.61\cdot 10^{\text{-}2}$ & $1.7\cdot 10^{\text{-}2}$ & $1.1\cdot 10^{\text{-}3}$ & $7.8\cdot 10^{\text{-}3}$ & -- & $9.0\cdot 10^{\text{-}3}$ & $5.3\cdot 10^{\text{-}3}$ & $3.0\cdot 10^{\text{-}4}$ & $5.0\cdot 10^{\text{-}7}$ \\
$1000$ & $20$ & $2.82\cdot 10^{\text{-}2}$ & $2.6\cdot 10^{\text{-}2}$ & $4.3\cdot 10^{\text{-}3}$ & $1.6\cdot 10^{\text{-}2}$ & -- & $1.1\cdot 10^{\text{-}3}$ & $7.4\cdot 10^{\text{-}8}$ & $3.0\cdot 10^{\text{-}4}$ & $8.2\cdot 10^{\text{-}8}$ \\
$2000$ & $20$ & $1.12\cdot 10^{\text{-}2}$ & $1.1\cdot 10^{\text{-}2}$ & $1.8\cdot 10^{\text{-}3}$ & $4.0\cdot 10^{\text{-}3}$ & -- & $8.2\cdot 10^{\text{-}3}$ & $6.7\cdot 10^{\text{-}7}$ & $1.5\cdot 10^{\text{-}6}$ & $2.3\cdot 10^{\text{-}10}$ \\
$1000$ & $25$ & $4.46\cdot 10^{\text{-}3}$ & $4.3\cdot 10^{\text{-}3}$ & $8.3\cdot 10^{\text{-}4}$ & $2.0\cdot 10^{\text{-}3}$ & -- & $3.5\cdot 10^{\text{-}3}$ & $7.3\cdot 10^{\text{-}14}$ & $3.8\cdot 10^{\text{-}7}$ & $5.1\cdot 10^{\text{-}7}$ \\
$2000$ & $25$ & $1.74\cdot 10^{\text{-}3}$ & $1.6\cdot 10^{\text{-}3}$ & $1.3\cdot 10^{\text{-}4}$ & $5.8\cdot 10^{\text{-}4}$ & -- & $1.4\cdot 10^{\text{-}3}$ & $9.1\cdot 10^{\text{-}12}$ & $1.6\cdot 10^{\text{-}8}$ & $1.3\cdot 10^{\text{-}11}$
\end{tabular}
\end{table}
\fi

\begin{table}[ht!]
\setlength{\tabcolsep}{3pt}
\centering
\caption{Comparison in terms of hypervolume for the knapsack problem.}
\label{table:ks_hv}
\begin{tabular}{r!{\vrule width 1.5pt}r!{\vrule width 1.5pt}r|r|r|r|r|r|r|r|r}
$T$ & $m$ & MONISE   & RAND   & NC       & PGEN     & RENN   & NSGA-II   & NSGA-III  & SPEA2    & SMS-E       \\
\noalign{\hrule height 1.5pt}
$1000$ & 5  & \cgs$5.38\nd 10^{\text{-}1}$ & \cgt$4.8\nd 10^{\text{-}1}$ & $3.6\nd 10^{\text{-}1}$ & $4.5\nd 10^{\text{-}1}$ & \cgf$5.39\nd 10^{\text{-}1}$ & $4.4\nd 10^{\text{-}1}$ & $4.5\nd 10^{\text{-}1}$ & $1.6\nd 10^{\text{-}2}$ & $1.6\nd 10^{\text{-}2}$ \\
$2000$ & 5  & \cgs$4.14\nd 10^{\text{-}1}$ & \cgf$4.2\nd 10^{\text{-}1}$ & $4.0\nd 10^{\text{-}1}$ & $3.9\nd 10^{\text{-}1}$ & \cgt$4.12\nd 10^{\text{-}1}$ & $2.7\nd 10^{\text{-}1}$ & $3.5\nd 10^{\text{-}1}$ & $4.7\nd 10^{\text{-}2}$ & $4.9\nd 10^{\text{-}2}$ \\
$1000$ & $10$ & \cgf$1.39\nd 10^{\text{-}1}$ & \cgs$1.2\nd 10^{\text{-}1}$ & $7.1\nd 10^{\text{-}2}$ & \cgt$1.1\nd 10^{\text{-}1}$ & $1.01\nd 10^{\text{-}1}$ & $9.6\nd 10^{\text{-}2}$ & $1.8\nd 10^{\text{-}2}$ & $9.7\nd 10^{\text{-}4}$ & $1.1\nd 10^{\text{-}3}$ \\
$2000$ & $10$ & \cgf$7.57\nd 10^{\text{-}2}$ & \cgs$6.9\nd 10^{\text{-}2}$ & $2.3\nd 10^{\text{-}2}$ & \cgt$6.6\nd 10^{\text{-}2}$ & $4.40\nd 10^{\text{-}2}$ & $3.5\nd 10^{\text{-}2}$ & $4.8\nd 10^{\text{-}3}$ & $5.5\nd 10^{\text{-}4}$ & $1.1\nd 10^{\text{-}4}$ \\
$1000$ & $15$ & \cgf$2.28\nd 10^{\text{-}2}$ & \cgs$1.8\nd 10^{\text{-}2}$ & $4.9\nd 10^{\text{-}3}$ & \cgt$1.3\nd 10^{\text{-}2}$ & -- & $1.2\nd 10^{\text{-}2}$ & $2.3\nd 10^{\text{-}3}$ & $1.2\nd 10^{\text{-}7}$ & $3.0\nd 10^{\text{-}7}$ \\
$2000$ & $15$ & \cgs$1.74\nd 10^{\text{-}2}$ & \cgf$1.8\nd 10^{\text{-}2}$ & $1.3\nd 10^{\text{-}3}$ & \cgt$8.8\nd 10^{\text{-}3}$ & -- & $6.0\nd 10^{\text{-}3}$ & $5.3\nd 10^{\text{-}3}$ & $5.9\nd 10^{\text{-}4}$ & $3.0\nd 10^{\text{-}6}$ \\
$1000$ & $20$ & \cgf$2.72\nd 10^{\text{-}2}$ & \cgs$2.5\nd 10^{\text{-}2}$ & $2.0\nd 10^{\text{-}3}$ & \cgt$1.5\nd 10^{\text{-}2}$ & -- & $1.1\nd 10^{\text{-}3}$ & $4.1\nd 10^{\text{-}8}$ & $2.8\nd 10^{\text{-}4}$ & $1.9\nd 10^{\text{-}5}$ \\
$2000$ & $20$ & \cgf$9.85\nd 10^{\text{-}3}$ & \cgs$9.6\nd 10^{\text{-}3}$ & $1.5\nd 10^{\text{-}3}$ & $3.5\nd 10^{\text{-}3}$ & -- & \cgf$6.5\nd 10^{\text{-}3}$ & $2.4\nd 10^{\text{-}8}$ & $6.2\nd 10^{\text{-}7}$ & $1.2\nd 10^{\text{-}6}$ \\
$1000$ & $25$ & \cgf$4.38\nd 10^{\text{-}3}$ & \cgs$4.2\nd 10^{\text{-}3}$ & $1.1\nd 10^{\text{-}3}$ & $2.0\nd 10^{\text{-}3}$ & -- & \cgt$3.4\nd 10^{\text{-}3}$ & $8.4\nd 10^{\text{-}14}$ & $5.0\nd 10^{\text{-}7}$ & $3.0\nd 10^{\text{-}9}$ \\
$2000$ & $25$ & \cgf$1.34\nd 10^{\text{-}3}$ & \cgs$1.3\nd 10^{\text{-}3}$ & $1.1\nd 10^{\text{-}4}$ & $4.9\nd 10^{\text{-}4}$ & -- & \cgt$1.1\nd 10^{\text{-}3}$ & $2.1\nd 10^{\text{-}12}$ & $1.6\nd 10^{\text{-}8}$ & $6.6\nd 10^{\text{-}14}$
\end{tabular}
\end{table}

\ifx
\begin{table}[ht!]
\setlength{\tabcolsep}{3pt}
\centering
\caption{Comparison in terms of execution time (in seconds) for the knapsack problem.}
\label{table:ks_time}
\begin{tabular}{r!{\vrule width 1.5pt}r!{\vrule width 1.5pt}r|r|r|r|r|r|r|r|r}
$T$ & $m$ & MONISE   & RAND   & NC       & PGEN     & RENN   & NSGA-II   & NSGA-III  & SPEA2    & SMS-E       \\
\noalign{\hrule height 1.5pt}
1000 & 5  & 11.4  & 0.2 & 4.6   & 1.6    & 33.4   & 3600.0 & 3600.1 & 3600.1 & 3600.0 \\
2000 & 5  & 2.1   & 0.1 & 1.6   & 1.0    & 37.7   & 3600.0 & 3600.0 & 3600.1 & 3600.0 \\
1000 & 10 & 57.2  & 0.5 & 8.8   & 3721.5 & 3646.9 & 3600.1 & 3601.7 & 3600.3 & 3600.0 \\
2000 & 10 & 18.5  & 0.4 & 10.5  & 3936.9 & 4646.0 & 3600.1 & 3601.0 & 3600.2 & 3600.0 \\
1000 & 15 & 112.5 & 0.9 & 93.7  & 5161.3 & --     & 3600.1 & 3610.5 & 3600.6 & 3601.5 \\
2000 & 15 & 69.3  & 0.7 & 107.3 & 5429.7 & --     & 3600.1 & 3607.0 & 3600.4 & 3603.2 \\
1000 & 20 & 173.5 & 2.0 & 200.7 & 6119.4 & --     & 3600.2 & 3664.3 & 3601.2 & 3600.1 \\
2000 & 20 & 125.6 & 1.2 & 118.9 & 5686.5 & --     & 3600.1 & 3646.7 & 3601.3 & 3602.5 \\
1000 & 25 & 235.9 & 2.3 & 443.7 & 7348.4 & --     & 3600.2 & 3747.1 & 3602.2 & 3617.3 \\
2000 & 25 & 173.4 & 2.1 & 258.4 & 7085.5 & --     & 3600.2 & 3852.7 & 3602.8 & 3956.9
\end{tabular}
\end{table}
\fi

\begin{table}[ht!]
\setlength{\tabcolsep}{3pt}
\centering
\caption{Comparison in terms of execution time (in seconds) for the knapsack problem.}
\label{table:ks_time}
\begin{tabular}{r!{\vrule width 1.5pt}r!{\vrule width 1.5pt}r|r|r|r|r|r|r|r|r}
$T$ & $m$ & MONISE   & RAND   & NC       & PGEN     & RENN   & NSGA-II   & NSGA-III  & SPEA2    & SMS-E       \\
\noalign{\hrule height 1.5pt}
1000 & 5  & \cgt2.3   & \cgf0.1 & 3.5   & \cgs1.2    & 25.2   & 3600.0 & 3600.0 & 3600.0 & 3600.0 \\
2000 & 5  & \cgt2.0   & \cgf0.1 & 1.6   & \cgs1.1    & 38.7   & 3600.0 & 3600.1 & 3600.1 & 3600.0 \\
1000 & 10 & \cgt17.6  & \cgf0.4 & \cgs8.8   & 3691.4 & 4017.7 & 3600.0 & 3601.4 & 3600.2 & 3600.0 \\
2000 & 10 & \cgt18.7  & \cgf0.4 & \cgs10.6  & 3802.1 & 3661.3 & 3600.0 & 3601.3 & 3600.2 & 3600.0 \\
1000 & 15 & \cgs67.1  & \cgf0.8 & \cgt96.4  & 4135.7 & --     & 3600.1 & 3613.3 & 3600.7 & 3600.0 \\
2000 & 15 & \cgs70.3  & \cgf0.7 & \cgt110.6 & 3878.1 & --     & 3600.1 & 3609.9 & 3600.4 & 3602.7 \\
1000 & 20 & \cgs120.5 & \cgf1.4 & \cgt178.9 & 6099.0 & --     & 3600.2 & 3663.7 & 3601.5 & 3600.1 \\
2000 & 20 & 125.3 & \cgf1.2 & \cgs118.0 & 5917.1 & --     & 3600.2 & 3630.7 & 3601.5 & 3619.0 \\
1000 & 25 & \cgs199.0 & \cgf2.2 & \cgt416.0 & 7448.5 & --     & 3600.2 & 3777.9 & 3602.4 & 3600.4 \\
2000 & 25 & \cgs175.4 & \cgf2.2 & \cgt267.5 & 7141.1 & --     & 3600.2 & 3852.7 & 3601.7 & 3628.0
\end{tabular}
\end{table}

\begin{table}[ht!]
\setlength{\tabcolsep}{3pt}
\centering
\caption{Comparison in terms of number of iterations (number of solutions for the deterministic algorithms and generations for the evolutionary) for the knapsack problem.}
\label{table:ks_iter}
\begin{tabular}{r!{\vrule width 1.5pt}r!{\vrule width 1.5pt}r|r|r|r|r|r|r|r|r}
$T$ & $m$ & MONISE   & RAND   & NC       & PGEN     & RENN   & NSGA-II   & NSGA-III  & SPEA2    & SMS-E       \\
\noalign{\hrule height 1.5pt}
1000 & 5  & 25  & 25  & 40  & 25 & 25 & 255369 & 100049 & 108247 & 1749084 \\
2000 & 5  & 25  & 25  & 35  & 25 & 25 & 256103 & 99892  & 117184 & 1888070 \\
1000 & 10 & 50  & 50  & 76  & 43 & 16 & 115841 & 3516   & 23500  & 217946  \\
2000 & 10 & 50  & 50  & 85  & 43 & 16 & 116506 & 3769   & 24304  & 163923  \\
1000 & 15 & 75  & 75  & 391 & 28 & 0  & 56233  & 349    & 9213   & 190051  \\
2000 & 15 & 75  & 75  & 415 & 28 & 0  & 63421  & 353    & 8638   & 8193    \\
1000 & 20 & 100 & 100 & 401 & 29 & 0  & 37626  & 50     & 4124   & 127563  \\
2000 & 20 & 100 & 100 & 262 & 29 & 0  & 40361  & 65     & 4333   & 7483    \\
1000 & 25 & 125 & 125 & 556 & 32 & 0  & 24199  & 16     & 2365   & 47205   \\
2000 & 25 & 125 & 125 & 322 & 32 & 0  & 24917  & 17     & 2172   & 2047   
\end{tabular}
\end{table}

\ifx
\begin{table}[ht!]
\setlength{\tabcolsep}{5pt}
\centering
\caption{Comparison in terms of hypervolume for the knapsack problem.}
\label{table:ks_hv}
\begin{tabular}{l|r|rrrrl}
\hline
             & $m$ & MONISE & NC    & NSGA-II & NSGA-III & SPEA2 \\ \hline
M\_05\_50\%  &  5  & $153\cdot 10^{-3}$    & $166\cdot 10^{-3}$    & $151\cdot 10^{-3}$    & $53\cdot 10^{-3}$    & $152\cdot 10^{-3}$    \\
M\_10\_50\% & 10 & $112\cdot 10^{-4}$   & $90\cdot 10^{-4}$   & $71\cdot 10^{-4}$   & $65\cdot 10^{-4}$   & $50\cdot 10^{-4}$   \\
M\_15\_50\% & 15 & $101\cdot 10^{-5}$  & $52\cdot 10^{-5}$  & $37\cdot 10^{-5}$  & $84\cdot 10^{-5}$  & $24\cdot 10^{-5}$  \\
M\_20\_50\% & 20 & $111\cdot 10^{-5}$ & $1\cdot 10^{-5}$ & $61\cdot 10^{-5}$ & $88\cdot 10^{-5}$ & $21\cdot 10^{-5}$   \\
\hline
\end{tabular}
\end{table}

\begin{table}[ht!]
\setlength{\tabcolsep}{5pt}
\centering
\caption{Comparison in terms of execution time (in seconds) for the knapsack problem.}
\label{table:ks_time}
\begin{tabular}{l|r|rrrrr}
\hline
            & $m$  & MONISE & NC    & NSGA-II & NSGA-III & SPEA2 \\ \hline
M\_05\_50\%  & 5  & 2   & 14    & 21  & 250    & 94   \\
M\_10\_50\% & 10 & 15  & 22    & 55  & 931    & 876   \\
M\_15\_50\% & 15 & 80  & 22    & 115 & 3,541  & 3,280 \\
M\_20\_50\% & 20 & 196 & 1,169 & 200 & 10,871 & 7,901  \\
\hline
\end{tabular}
\end{table}
\fi

\subsection{Comparative Analysis}

Except for the RAND approach, the execution time for MONISE is lower than those produced by the other methods in the majority of the instances. Basically, the experimental design offered more computational time to the other methods in an attempt to reach or to slightly overcome, when possible, the performance of MONISE. The time limit of 3600 seconds was exploited by almost all competing methods. It justifies the choice of a relatively low number of desired solutions ($5\times m$), since any increase in this number would imply in an even lower hypervolume performance for those methods.

It is noticeable that MONISE exhibits a consistent performance on both problems and metrics. The performance of MONISE with respect to the hypervolume metric seems to be insensitive to the number of objectives, and, in the convex nonlinear problem (multilabel classification) MONISE is the best method in 3, and second best in 2 out of 5 instances. For the combinatorial problem (knapsack), MONISE  is the best method in 7, and the second best in 3 out of 10 instances. When MONISE is not the best method, the hypervolume is quite close to the best, with a comparatively reduced execution time, except for the RAND approach. The good performance of the algorithms based on weighted sum method scalarization for the knapsack problem can be explained by the relatively low number of explored solutions ($5\times m$). Indeed, a rough exploration of solutions attainable by this scalarization already leads to a good performance in the hypervolume.

Analysing the competing methods, it is possible to see that a ``smart'' random method (RAND) is capable of achieving good results, coming closer to MONISE when the solutions are scattered in the knapsack problem, for example. However, improved results are produced by more informed methods (MONISE and RENN). As expected, the only scenario where RAND achieves the best results is in the knapsack problem with loose enough constraints ($T=2000$) responsible for increasing the randomness of the multi-objective problem. The NC method seems to be more appropriate for lower dimensions, exhibiting hypervolume and execution time issues in higher dimensions. PGEN is not scalable but seems to be a reasonable method, although loses to RAND in all situations. This behaviour seems to be caused by its heuristic approach when a weight vector with negative values is selected. Additionally, it can be noticed that, with the increase in the number of objectives, the number of solutions found is reduced due to high cost of enumerating the convex-hull, which precludes the search of new solutions within the designated time (3600 seconds). RENN has a high performance in low dimensional problems at the price of a quite high computational cost, preventing this method from being applicable in moderate to high-dimensional scenarios. This cost is caused by the high dimension of the convex hull associated with the high number of (dummy-)points. As well as in PGEN, the high cost of enumerating the convex-hull also precludes the search of new solutions within the designated time (3600 seconds) in some situations. The last two polyhedral methods have a high cost and scalability issues caused by the enumerative nature of convex hull calculation.

The evolutionary algorithms were never capable of achieving the best performance, even taking 3600 seconds. For the multilabel classification problem, it is noticeable that NSGA-II is the best method for lower dimensions and NSGA-III for higher dimensions (except for genbase). For the knapsack problem, all methods but NSGA-II have their performance reduced for higher dimensions. This better performance of NSGA-II might be explained by its lower computational cost per iteration, allowing NSGA-II to make more iterations than other evolutionary algorithms for the same limit of computational time. Notice that SMS-E only evaluates a single individual at each iteration, and yet NSGA-II has a similar number of iterations (or more) on high dimensions.

Observing the results of the ``smart'' random weight vector sampling (RAND) we can extract three relevant characteristics of the weighted sum method: (1) even with simple methods of weight vector sampling, this scalarization can have a good performance in the hypervolume metric; (2) the computational cost is quite reduced; and (3) the computational cost is well behaved with the increase in the number of objectives. Exploring these benefits, MONISE spends additional computational resources to improve the search for specific weight vectors, thus improving the hypervolume indicator.

The results show that MONISE was capable of maintaining the characteristics of the weighted sum method with random sampling, but improving the hypervolume indicator at the cost of a slight increase in the computational cost. The time limitation in the weight vector calculation controls the behavior of the execution time in any scale. Despite this time limitation, the comparative performance in terms of hypervolume was always consistent, not being the best method in only 3 out of 15 instances, keeping a good hypervolume performance in those cases and always having an outstanding performance in terms of execution time. The reduced computational time required to calculate the weight vectors should be attributed to mixed-integer formulation and its well-consolidated methodologies. Very efficient commercial softwares (as well as some open source softwares) are capable of achieving good solutions even with small execution time.

The overall experimental behavior of MONISE is outstanding. By avoiding the convex hull calculation, it consistently avoids the high computational cost required by PGEN and RENN, which took all the available computational time in some cases. Also, the MONISE strategy to avoid and correct negative weight vectors led to better quality solutions than PGEN -- with lower computational cost in the majority of the cases -- and led to lower computational cost when compared to RENN -- a method whose executions were broken in high dimensions due to excess of memory usage. These results give support to endorse the weighted sum method, associated with MONISE, as a useful tool for high-dimensional multi-objective problems.

\section{Conclusion} \label{sec:conclusion}

The main goal of this work was to extend the Non-Inferior Set Estimation (NISE) algorithm for more than two objectives. The new proposal was called Many Objective NISE or simply MONISE. The methodology is based on the well known scalarization known as the weighted sum method. On one hand, this aspect contributes to the generality, simplicity and interpretability of the work, which can be briefly described as a computationally efficient procedure devoted to a recursive and automatic sampling of the Pareto frontier. On the other hand, in the case of non-convex Pareto frontiers, efficient solutions dominated by convex combinations of other efficient solutions are not achievable by MONISE, as it happens for all methods based on the weighted sum method.

Relevant theoretical aspects have been demonstrated for the NISE algorithm considering two-objectives problems, due to the convergence behavior of the weighted sum method. However, using again the behavior of the weighted sum method, it was demonstrated that those properties are no longer valid for three or more objectives, using a three-objective optimization problem instance as a case study (see Appendix \ref{appendix:nise_nonrecursivity}).

Theoretical developments, together with the proper characterization of the procedure of the inner / outer approximation approach adopted by the weighted sum method to estimate the Pareto frontier, allowed a deeper understanding of the convergence behavior of the weighted sum method, as well as of the NISE method. Those insights act as the ground motivation to propose MONISE as an extension of NISE capable of dealing with two or more objectives. In practical terms, the original procedure in NISE responsible for recursively finding the weighting vectors was replaced by a mixed-integer linear programming formulation (see Definition \ref{def:moo:monise_w_estim_pli}) for which efficient commercial and open-source solvers are available. Under a geometrical perspective, MONISE operates by iteratively finding weighting vectors associated with the highest difference between the inner and outer approximation frontiers.

The empirical results included comparisons with four heuristically-based methods, one random approach and three non-heuristic multi-objective optimization algorithms, which reveals a consistent performance for the MONISE proposal. Despite not being the best-evaluated method all the time in terms of hypervolume, the deterministic and recursive procedure used to sample the Pareto frontier seems to confer robustness and consistency to the proposal. MONISE also scales favorably with the increase in the number of objectives, being always efficient in terms of computational cost. We could demonstrate the superiority of MONISE over methods that enumerate the convex hull facets. The computational overload caused by the enumeration is mitigated by our mixed integer modeling of the multi-objective search strategy, associated with an early stopping of the search that finds a high-quality facet without the enumeration of all the convex hull facets. Notice that we are not necessarily solving a multi-objective integer and mixed-integer program. Instead, we are simply resorting to a mixed-integer linear solver to cleverly sample the Pareto frontier of a great variety of multi-objective optimization problems.

Further developments in this area reside basically in two main fronts: (1) adaptation of the algorithm to take advantage of scalarization and development of dedicated solvers for the single-objective optimization problems associated with each weighting vector. For instance, we may try to estimate the lower bound of a sub-optimal solution, so that this lower bound can be used in outer approximations; (2) refinement of the mechanisms involving the weighted sum method. One possibility is to make a better usage of the calls to the scalarization inside the mixed-integer linear solver. The idea is to insert the related constraints as soon as a promising weighting vector is found by the mixed-integer linear solver.

\bibliographystyle{elsarticle}
\biboptions{authoryear}

\section*{Acknowledgements}

This work was supported by grants from CNPq (process \#307228/2018-5) and FAPESP\--CAPES (process \#2014/13533-0).

\bibliography{main}

\begin{thebibliography}{47}
\expandafter\ifx\csname natexlab\endcsname\relax\def\natexlab#1{#1}\fi
\expandafter\ifx\csname url\endcsname\relax
  \def\url#1{\texttt{#1}}\fi
\expandafter\ifx\csname urlprefix\endcsname\relax\def\urlprefix{URL }\fi

\bibitem[{Bazgan et~al.(2009)Bazgan, Hugot, and Vanderpooten}]{Bazgan2009a}
Bazgan, C., Hugot, H., Vanderpooten, D., 2009. {Solving efficiently the 0-1 multi-objective knapsack problem}. Computers {\&} Operations Research 36~(1), 260--279.

\bibitem[{Benson(1998)}]{Benson1998}
Benson, H.~P., 1998. {An Outer Approximation Algorithm for Generating All Efficient Extreme Points in the Outcome Set of a Multiple Objective Linear Programming Problem}. Journal of Global Optimization 13~(1), 1--24.

\bibitem[{Beume et~al.(2007)Beume, Naujoks, and Emmerich}]{Beume2007}
Beume, N., Naujoks, B., Emmerich, M., 2007. {SMS-EMOA: Multiobjective selection based on dominated hypervolume}. European Journal of Operational Research 181~(3), 1653--1669.

\bibitem[{Bokrantz and Forsgren(2013)}]{Bokrantz2013}
Bokrantz, R., Forsgren, A., 2013. {An algorithm for approximating convex pareto surfaces based on dual techniques}. INFORMS Journal on Computing 25~(2), 377--393.

\bibitem[{Burachik et~al.(2013)Burachik, Kaya, and Rizvi}]{Burachik2013}
Burachik, R.~S., Kaya, C.~Y., Rizvi, M.~M., 2013. {A New Scalarization Technique to Approximate Pareto Fronts of Problems with Disconnected Feasible Sets}. Journal of Optimization Theory and Applications 162, 428--446.

\bibitem[{Caballero and Hern{\'{a}}ndez(2004)}]{Caballero2004}
Caballero, R., Hern{\'{a}}ndez, M., 2004. {The controlled estimation method in the multiobjective linear fractional problem}. Computers {\&} Operations Research 31~(11), 1821--1832.

\bibitem[{Ceyhan et~al.(2019)Ceyhan, K{\"{o}}ksalan, and Lokman}]{Ceyhan2019}
Ceyhan, G., K{\"{o}}ksalan, M., Lokman, B., 2019. {Finding a representative nondominated set for multi-objective mixed integer programs}. European Journal of Operational Research 272~(1), 61--77.

\bibitem[{Cohon et~al.(1979)Cohon, Church, and Sheer}]{Cohon1979}
Cohon, J.~L., Church, R.~L., Sheer, D.~P., 1979. {Generating multiobjective trade‐offs: An algorithm for bicriterion problems}. Water Resources Research 15~(5), 1001--1010.

\bibitem[{Craft et~al.(2006)Craft, Halabi, Shih, and Bortfeld}]{Craft2006}
Craft, D.~L., Halabi, T.~F., Shih, H.~A., Bortfeld, T.~R., 2006. {Approximating convex Pareto surfaces in multiobjective radiotherapy planning}. Medical Physics 33~(9), 3399--3407.

\bibitem[{Das and Dennis(1997)}]{Das1997}
Das, I., Dennis, J., 1997. {A closer look at drawbacks of minimizing weighted sums of objectives for Pareto set generation in multicriteria optimization problems}. Structural Optimization 14~(1), 63--69.

\bibitem[{Das and Dennis(1998)}]{Das1998}
Das, I., Dennis, J., 1998. {Normal-boundary intersection: A new method for generating the Pareto surface in nonlinear multicriteria optimization problems}. SIAM Journal on Optimization 8~(3), 631--657.

\bibitem[{Deb and Jain(2013)}]{Deb2013}
Deb, K., Jain, H., 2013. {An Evolutionary Many-Objective Optimization Algorithm Using Reference-point Based Non-dominated Sorting Approach, Part I: Solving Problems with Box Constraints}. IEEE Transactions on Evolutionary Computation 18~(4), 577--601.

\bibitem[{Deb et~al.(2002)Deb, Pratap, Agarwal, and Meyarivan}]{Deb2002}
Deb, K., Pratap, A., Agarwal, S., Meyarivan, T., 2002. {A fast and elitist multiobjective genetic algorithm: NSGA-II}. IEEE Transactions on Evolutionary Computation 6~(2), 182--197.

\bibitem[{Ehrgott et~al.(2011)Ehrgott, Shao, and Sch{\"{o}}bel}]{Ehrgott2011}
Ehrgott, M.~M., Shao, L.~L., Sch{\"{o}}bel, A.~A., 2011. {An approximation algorithm for convex multi-objective programming problems}. Journal of Global Optimization 50~(3), 397--416.

\bibitem[{Eichfelder(2009{\natexlab{a}})}]{Eichfelder2009}
Eichfelder, G., 2009{\natexlab{a}}. {An Adaptive Scalarization Method in Multiobjective Optimization}.

\bibitem[{Eichfelder(2009{\natexlab{b}})}]{Eichfelder2009a}
Eichfelder, G., 2009{\natexlab{b}}. {Scalarizations for adaptively solving multi-objective optimization problems}. Computational Optimization and Applications 44~(2), 249--273.

\bibitem[{Fleischer(2003)}]{Fleischer2003}
Fleischer, M., 2003. {The measure of Pareto optima applications to multi-objective metaheuristics}. Evolutionary Multi-Criterion Optimization 17, 519--533.

\bibitem[{Geoffrion(1968)}]{Geoffrion1968}
Geoffrion, A.~M., 1968. {Proper efficiency and the theory of vector maximization}. Journal of Mathematical Analysis and Applications 22~(3), 618--630.

\bibitem[{Ishibuchi et~al.(2009)Ishibuchi, Sakane, Tsukamoto, and Nojima}]{Ishibuchi2009}
Ishibuchi, H., Sakane, Y., Tsukamoto, N., Nojima, Y., 2009. {Evolutionary many-objective optimization by NSGA-II and MOEA/D with large populations}. Conference Proceedings - IEEE International Conference on Systems, Man and Cybernetics 1, 1758--1763.

\bibitem[{Khorram et~al.(2014)Khorram, Khaledian, and Khaledyan}]{Khorram2014}
Khorram, E., Khaledian, K., Khaledyan, M., 2014. {A numerical method for constructing the Pareto front of multi-objective optimization problems}. Journal of Computational and Applied Mathematics 261, 158--171.

\bibitem[{Kim et~al.(2006)Kim, Weck, and De~Weck}]{Kim2005}
Kim, I., Weck, O., De~Weck, O., 2006. {Adaptive weighted sum method for multiobjective optimization: a new method for Pareto front generation}. Structural and Multidisciplinary Optimization 31~(2), 105--116.

\bibitem[{Kirlik and Sayin(2014)}]{Kirlik2014c}
Kirlik, G., Sayin, S., 2014. {A new algorithm for generating all nondominated solutions of multiobjective discrete optimization problems}. European Journal of Operational Research 232~(3), 479--488.

\bibitem[{Klamroth et~al.(2003)Klamroth, Tind, and Wiecek}]{Klamroth2003}
Klamroth, K., Tind, J., Wiecek, M.~M., 2003. {Unbiased approximation in multicriteria optimization}. Mathematical Methods of Operations Research 56~(3), 413--437.

\bibitem[{Koski(1985)}]{Koski1985}
Koski, J., 1985. {Defectiveness of weighting method in multicriterion optimization of structures}. Communications in Applied Numerical Methods 1~(May), 333--337.

\bibitem[{Kukkonen et~al.(2007)Kukkonen, Member, and Lampinen}]{Kukkonen2007}
Kukkonen, S., Member, S., Lampinen, J., 2007. {Ranking-Dominance and Many-Objective Optimization}. In: IEEE Congress on Evolutionary Computation. pp. 3983--3990.

\bibitem[{Marler and Arora(2004)}]{Marler2004}
Marler, R., Arora, J., 2004. {Survey of multi-objective optimization methods for engineering}. Structural and Multidisciplinary Optimization 26~(6), 369--395.

\bibitem[{Marler and Arora(2010)}]{Marler2009}
Marler, R., Arora, J., 12 2010. {The weighted sum method for multi-objective optimization: new insights}. Structural and Multidisciplinary Optimization 41~(6), 853--862.

\bibitem[{Masin and Bukchin(2008)}]{Masin2008a}
Masin, M., Bukchin, Y., 2008. {Diversity maximization approach for multiobjective optimization}. Operations Research 56~(2), 411--424.

\bibitem[{Messac et~al.(2003)Messac, Ismail-Yahaya, and Mattson}]{Messac2003}
Messac, A., Ismail-Yahaya, A., Mattson, C., 7 2003. {The normalized normal constraint method for generating the Pareto frontier}. Structural and Multidisciplinary Optimization 25~(2), 86--98.

\bibitem[{Messac and Mattson(2004)}]{Messac2004}
Messac, A., Mattson, C., 2004. {Normal constraint method with guarantee of even representation of complete Pareto frontier}. AIAA Journal 42~(10), 2101--2111.

\bibitem[{Michalewicz and Arabas(1994)}]{Michalewicz1994}
Michalewicz, Z., Arabas, J., 1994. {Genetic algorithms for the 0/1 knapsack problem}. Methodologies for Intelligent Systems 869, 134--143.

\bibitem[{Miettinen(1999)}]{Miettinen1999}
Miettinen, K., 1999. {Nonlinear Multiobjective Optimization}. Springer.

\bibitem[{Nobakhtian and Shafiei(2016)}]{Nobakhtian2017}
Nobakhtian, S., Shafiei, N., 2016. {A Benson type algorithm for nonconvex multiobjective programming problems}. Top 25~(2), 271--287.

\bibitem[{{\"{O}}zlen and Azizoglu(2009)}]{Ozlen2009}
{\"{O}}zlen, M., Azizoglu, M.~M., 2009. {Multi-objective integer programming: A general approach for generating all non-dominated solutions}. European Journal of Operational Research 199~(1), 25--35.

\bibitem[{{\"{O}}zlen et~al.(2014){\"{O}}zlen, Burton, and MacRae}]{Ozlen2014}
{\"{O}}zlen, M., Burton, B., MacRae, C., 2014. {Multi-Objective Integer Programming: An Improved Recursive Algorithm}. Journal of Optimization Theory and Applications 160, 470--482.

\bibitem[{{\"{O}}zpeynirci and K{\"{o}}ksalan(2010)}]{Ozpeynirci2010a}
{\"{O}}zpeynirci, {\"{O}}., K{\"{o}}ksalan, M., 2010. {An exact algorithm for finding extreme supported nondominated points of multiobjective mixed integer programs}. Management Science 56~(12), 2302--2315.

\bibitem[{Rennen et~al.(2009)Rennen, van Dam, and den Hertog}]{Rennen2009}
Rennen, G., van Dam, E., den Hertog, D., 2009. {Enhancement of Sandwich Algorithms for Approximating Higher Dimensional Convex Pareto Sets}. Ssrn~(December 2018).

\bibitem[{Romero and Rehman(2003)}]{Romero2003}
Romero, C., Rehman, T., 2003. {Multiobjective programming}. In: Romero, C., Rehman, T. (Eds.), Multiple Criteria Analysis for Agricultural Decisions. Vol.~11 of Developments in Agricultural Economics. Elsevier, Ch.~4, pp. 47--61.

\bibitem[{Ryu et~al.(2009)Ryu, Kim, and Wan}]{Ryu2009a}
Ryu, J.~H., Kim, S., Wan, H., 2009. {Pareto front approximation with adaptive weighted sum method in multiobjective simulation optimization}. Proceedings - Winter Simulation Conference, 623--633.

\bibitem[{Sanchis et~al.(2007)Sanchis, Mart{\'{i}}nez, Blasco, and Salcedo}]{Sanchis2007}
Sanchis, J., Mart{\'{i}}nez, M., Blasco, X., Salcedo, J.~V., 2007. {A new perspective on multiobjective optimization by enhanced normalized normal constraint method}. Structural and Multidisciplinary Optimization 36~(5), 537--546.

\bibitem[{Shao and Ehrgott(2008)}]{Shao2008}
Shao, L., Ehrgott, M., 2008. {Approximately solving multiobjective linear programmes in objective space and an application in radiotherapy treatment planning}. Mathematical Methods of Operations Research 68~(2), 257--276.

\bibitem[{Smith and Tromble(2004)}]{Smith2004}
Smith, N.~A., Tromble, R.~W., 2004. {Sampling Uniformly from the Unit Simplex Naive Algorithms}. Johns Hopkins University Tech. Rep~(29), 1--6.

\bibitem[{Snyder and ReVelle(1997)}]{Snyder1997}
Snyder, S., ReVelle, C., 1997. {Multiobjective grid packing model: an application in forest management}. Location Science 5~(3), 165--180.

\bibitem[{Solanki et~al.(1993)Solanki, Appino, and Cohon}]{Solanki1993}
Solanki, R.~S., Appino, P.~A., Cohon, J.~L., 1993. {Approximating the noninferior set in multiobjective linear programming problems}. European Journal of Operational Research 68~(3), 356--373.

\bibitem[{Sylva and Crema(2004)}]{Sylva2004}
Sylva, J., Crema, A., 2004. {A method for finding the set of non-dominated vectors for multiple objective integer linear programs}. European Journal of Operational Research 158~(1), 46--55.

\bibitem[{Xia et~al.(2018)Xia, Vera, and Zuluaga}]{Xia2015}
Xia, W., Vera, J., Zuluaga, L.~F., 2018. {Globally solving Non-Convex Quadratic Programs via Linear Integer Programming techniques}. arXiv preprint arXiv:1511.02423, 1--17.

\bibitem[{Zitzler et~al.(2001)Zitzler, Laumanns, and Thiele}]{Zitzler2001}
Zitzler, E., Laumanns, M., Thiele, L., 2001. {SPEA2: Improving the Strength Pareto Evolutionary Algorithm}. Evolutionary Methods for Design Optimization and Control with Applications to Industrial Problems, 95--100.

\end{thebibliography}

\clearpage
\setcounter{page}{1}
\begin{center}
\Large Supplementary Material:\\ An Extension of the Non-Inferior Set Estimation Algorithm\\ for Many Objectives
\end{center}

\begin{appendices}
\section{Recursivity of the NISE algorithm and its associated properties} \label{appendix:nise_recursivity}

Due to the fact that the NISE algorithm works locally in the objective space to obtain new solutions from already obtained solutions, some properties are necessary to ensure recursivity. Solutions in the objective space will be denoted by $\vec{r}$ (and not $\vec{f}(\vec{x}))$ due to notation simplicity and also to allow the representation of efficient solutions in the objective space that may not correspond to feasible solutions in the decision space.
Consider a neighborhood formed by the solutions $\vec{r}^1$ and $\vec{r}^2$, obtained using the weighting vectors $\vec{w}^1$ and $\vec{w}^2$. To ensure the recursivity of the algorithm, the segment that contains $\vec{r}^1$ and $\vec{r}^2$, having $\vec{w}$ as normal vector, should always generate a solution between $\vec{r}^1$ and $\vec{r}^2$.

Firstly, we show that a new weighting vector $\vec{w}$ must be a convex combination of $\vec{w}^1$ and $\vec{w}^2$, and then we show that a linear combination of $\vec{w}^1$ and $\vec{w}^2$ generates a solution between $\vec{r}^1$ and $\vec{r}^2$.

\begin{theorem}\label{theo:moo:weight_2Drecur} \textbf{Recursivity of the weighted sum method for two dimensions}

Given two efficient solutions $\vec{r}^{1}$ and $\vec{r}^{2}$ and the corresponding weighting vectors used to find them, $\vec{w}^{1}$ and $\vec{w}^{2}$, where $w^1_i,w^2_i$ $>0\ \forall i \in \{1,2\}$ and $w^1_1+w^1_2=w^2_1+w^2_2=1$, the optimality of the weighted sum method implies:

\begin{subequations}
\begin{align}
{\vec{w}^{1}}^\top\vec{r}^{1} < {\vec{w}^{1}}^\top\vec{r}, \forall\ \vec{r} \in \Psi, \label{theo:moo:sens-c}\\
{\vec{w}^{2}}^\top\vec{r}^{2} < {\vec{w}^{2}}^\top\vec{r}, \forall\ \vec{r} \in \Psi. \label{theo:moo:sens-d}
\end{align}
\end{subequations}

Then there exist $\alpha,\beta>0, \alpha+\beta=1$ such that:

\begin{equation} \label{eq:thet:aux1}
\begin{cases}
\vec{w}^\top \vec{y}^1 = \vec{w}^\top \vec{y}^2\\
\vec{w} = \alpha \vec{w}^1 + \beta \vec{w}^2\\
w_1+w_2 = 1.
\end{cases}
\end{equation}

\begin{proof} The proof by contradiction will be subdivided into three cases:
\begin{description}
\item[Case 1:] $\beta\not = (1-\alpha)$.

Since $\vec{w}\in \mathbb{R}^2$ and $\vec{w}^1\not = \vec{w}^2$, $\{\vec{w}^1,\vec{w}^2\}$ is basis of $\mathbb{R}^2$, then $\exists \alpha, \beta: \vec{w}=\alpha\vec{w}^1+\beta\vec{w}^2$. Notice that $\vec{w}_1+\vec{w}_2 = \alpha(w^1_1+w^1_2)+\beta(w^2_1+w^2_2) = \alpha+\beta$. However, if $\beta\not = (1-\alpha)$, $\vec{w}_1+\vec{w}_2 \not = 1$, thus generating a contradiction.

\item[Case 2:] $\beta < 0$.

Multiplying (\ref{theo:moo:sens-c}) by $\alpha$ and (\ref{theo:moo:sens-d}) by $\beta$:

\begin{equation*}
\begin{cases}
\alpha{\vec{w}^{1}}^\top\vec{r}^{1} < \alpha{\vec{w}^{1}}^\top\vec{r}, \forall\ \vec{r} \in \Psi\\
\beta {\vec{w}^{2}}^\top\vec{r} < \beta {\vec{w}^{2}}^\top\vec{r}^{2}, \forall\ \vec{r} \in \Psi.
\end{cases}
\end{equation*}

We obtain $\vec{w}^\top \vec{r}^1 < \vec{w}^\top \vec{r}^2$, which generates a contradiction.

\item[Case 3:] $\alpha <0$ - Multiplying (\ref{theo:moo:sens-c}) by $\alpha$ and (\ref{theo:moo:sens-d}) by $\beta$:

\begin{equation*}
\begin{cases}
\alpha{\vec{w}^{1}}^\top\vec{r}^{1} > \alpha{\vec{w}^{1}}^\top\vec{r}, \forall\ \vec{r} \in \Psi\\
\beta {\vec{w}^{2}}^\top\vec{r} > \beta {\vec{w}^{2}}^\top\vec{r}^{2}, \forall\ \vec{r} \in \Psi.
\end{cases}
\end{equation*}

We obtain $\vec{w}^\top \vec{r}^1 > \vec{w}^\top \vec{r}^2$, which generates a contradiction.

\end{description}
Since all cases generate a contradiction, the proof is concluded.
\end{proof}

\end{theorem}

To show that a linear combination of $\vec{w}^1$ and $\vec{w}^2$ generates a solution between its associated solutions $\vec{r}^1$ and $\vec{r}^2$, it is necessary to prove that, when we increase the weighting of an objective function (in a two dimensional space) its objective value decreases.

\begin{theorem}\label{theo:moo:weight_2Dsens} \textbf{Monotonicity of the weighted sum method for two dimensions}

Given two efficient solutions $\vec{r}^{1}$ and $\vec{r}^{2}$ and their corresponting weighting vectors $\vec{w}^{1}$ and $\vec{w}^{2}$ ($\vec{w}^{1} \not = \vec{w}^{2}$), where $w^1_i,w^2_i >0\ \forall i \in \{1,2\}$ and $w^1_1+w^1_2=w^2_1+w^2_2=1$, the optimality of the weighted sum method implied that:

\begin{subequations}
\begin{align}
{\vec{w}^{1}}^\top\vec{r}^{1} \leq {\vec{w}^{1}}^\top\vec{r}^2, \label{theo:moo:sens-a}\\
{\vec{w}^{2}}^\top\vec{r}^{2} \leq {\vec{w}^{2}}^\top\vec{r}^1. \label{theo:moo:sens-b}
\end{align}
\end{subequations}

If $w^1_1>w^2_1$ ($w^1_1<w^2_1$) then $r^1_1 \leq r^2_1$ ($r^1_1 \geq r^2_1$) and $r^1_2 \geq r^2_2$ ($r^1_2 \leq r^2_2$).

\begin{proof}
The proof is again by contradiction:
\begin{description}
\item[Case 1:] $r^1_1>r^2_1$ and $r^1_2\geq r^2_2$: 

Expanding (\ref{theo:moo:sens-a}), we have:

\begin{align*}
w^{1}_1\left(r^{1}_1-r^{2}_1\right) \leq w^{1}_2\left(r^{2}_2-r^{1}_2\right),\\
\frac{w^{1}_1}{w^{1}_2} \leq \frac{r^{2}_2-r^{1}_2}{r^{1}_1-r^{2}_1} \leq 0.
\end{align*}

However $\frac{w^{1}_1}{w^{1}_2}>0$ due to the positivity of the $w_i^j$s, leading to a contradiction.

\item[Case 2:] $r^1_1\leq r^2_1$ and $r^1_2<r^2_2$:

Expanding (\ref{theo:moo:sens-b}), we have:
\begin{align*}
w^{2}_2\left(r^{2}_2-r^{1}_2\right) \leq w^{2}_1\left(r^{1}_1-r^{2}_1\right),\\
\frac{w^{2}_1}{w^{2}_1} \leq \frac{r^{1}_1-r^{2}_1}{r^{2}_2-r^{1}_2} \leq 0.
\end{align*}

However $\frac{w^{2}_2}{w^{2}_1}>0$ due to the positivity of the $w_i^j$s, leading to a contradiction.

\item[Case 3:] $r^1_1>r^2_1$ and $r^1_2<r^2_2$:

Expanding (\ref{theo:moo:sens-a}), we obtain:
\begin{equation*}
\frac{w^{1}_1}{w^{1}_2} \leq \frac{r^{2}_2-r^{1}_2}{r^{1}_1-r^{2}_1}.
\end{equation*}

Expanding (\ref{theo:moo:sens-b}), we get:
\begin{equation*}
\frac{r^{2}_2-r^{1}_2}{r^{1}_1-r^{2}_1} \leq \frac{w^{2}_1}{w^{2}_2}.
\end{equation*}

Therefore we get $\frac{w^{2}_1}{w^{2}_2} \geq \frac{w^{1}_1}{w^{1}_2} \geq 1$, implying that $w^{2}_1 \geq w^{2}_2$, leading to a contradiction.
\end{description}

Since all cases generate a contradiction, the proof is concluded.
\end{proof}

\end{theorem}

\begin{theorem}\label{theo:moo:weight_2Dlocal} \textbf{Locality of the weighted sum method for two dimensions}

Given two efficient solutions $\vec{r}^{1}$ and $\vec{r}^{2}$ and their respective weighting vectors $\vec{w}^{1}$ and $\vec{w}^{2}$, where $w^1_i,w^2_i >0\ \forall i \in \{1,2\}$ and $w^1_1+w^1_2=w^2_1+w^2_2=1$, the optimality of the weighted sum method implies that:

\begin{subequations}
\begin{align*}
{\vec{w}^{1}}^\top\vec{r}^{1} \leq {\vec{w}^{1}}^\top\vec{r}, \forall\ \vec{r} \in \Omega,\\
{\vec{w}^{2}}^\top\vec{r}^{2} \leq {\vec{w}^{2}}^\top\vec{r}, \forall\ \vec{r} \in \Omega.
\end{align*}
\end{subequations}

If a solution $\overline{\vec{r}}$ is generated by the weighting vector $\vec{w} = (1-\alpha)\vec{w}^1 + \alpha \vec{w}^2$, $0 < \alpha < 1$, then:

\begin{equation*}
\min(r^1_i,r^2_i) \leq \overline{\vec{r}}_i \leq \max(r^1_i,r^2_i),\ i \in \{1,2\}.
\end{equation*}

\begin{proof}
Without loss of generality, it is supposed that $w^2_1 >w_1 > w^1_1$. Since we have that $w_1 > w^1_1$ and $w_2 < w^1_2$, using Theorem \ref{theo:moo:weight_2Dsens} we conclude that:

\begin{align*}
\overline{r}_1 \leq r^1_1, \\
\overline{r}_2 \geq r^1_2.
\end{align*}

Furthermore, we have that $w_1 < w^2_1$ and $w_2 > w^2_2$, and using Theorem \ref{theo:moo:weight_2Dsens} we conclude that:

\begin{align*}
r^2_1 \leq \overline{r}_1,\\
r^2_2 \geq \overline{r}_2.
\end{align*}

Therefore we conclude that $\min(r^1_1,r^2_1) \leq \overline{r}_1 \leq \max(r^1_1,r^2_1)$ and $\min(r^1_2,r^2_2) \leq \overline{r}_2 \leq \max(r^1_2,r^2_2)$.
\end{proof}

\end{theorem}

Supported by these properties of the weighted sum method, for any two efficient solutions with their associated weighting vectors, the weighting vector that solves (\ref{eq:moo:nise_w_estim}) is a convex combination of the previous weighting vectors (Theorem \ref{theo:moo:weight_2Drecur}), and the solution found for this new weighting vector remains between the two previous solutions (Theorem \ref{theo:moo:weight_2Dlocal}). These properties lead to a convergent procedure for the NISE algorithm in two dimensions. However, when the number of objective functions is superior to two, those properties are not necessarily preserved, as demonstrated in Appendix \ref{appendix:nise_nonrecursivity}.

\section{Violation of relevant properties of NISE in three or more dimensions} \label{appendix:nise_nonrecursivity}

A convex multi-objective problem with three objectives will be used to illustrate the violation of relevant properties of NISE when the number of dimensions is superior to two. Counter-examples for some conditions used to support the proofs in Appendix \ref{appendix:nise_recursivity} are presented. Consider the multi-objective problem:

\begin{equation} \label{theo:moo:cex3D}
\begin{aligned}
& \underset{\vec{x}}{\text{minimize}}
& & \vec{r} = \vec{f}(\vec{x}) = \left[x_1^2,x_2^2,x_3^2\right]\\
& \text{subject to}
& & \vec{x}^\top \vec{1} = 1,\ x_1, x_2, x_3 \geq 0.
\end{aligned}
\end{equation}

\begin{theorem}\label{theo:moo:weight_3Dnonrecur} \textbf{Non-recursivity of the weighted sum method in dimensions superior to two.}

Given the problem in (\ref{theo:moo:cex3D}), when the weighting vectors  $\vec{w}^1 = (0.24, 0.68, 0.08)$, $\vec{w}^2 = (0.23, 0.5, 0.27)$ and $\vec{w}^3 = (0.17, 0.38, 0.45)$ are used in the weighted sum method, the following solutions are obtained: $\vec{r}^1 = (0.05$, $0.006, 0.47)$, $\vec{r}^2 = (0.18, 0.04, 0.13)$, $\vec{r}^3 = (0.3, 0.06, 0.04)$. A hyperplane defined by these solutions is given by $\vec{w}^\top \vec{y} = (-0.14, 1.08, 0.06)\vec{y} = 0.2545,\ \vec{y} \in \mathbb{R}^3$, where the first component of the normal vector is negative $w_1 = -0.14$, thus breaking the recursivity of the NISE method for a dimension superior to two.
\end{theorem}

\begin{theorem}\label{theo:moo:weight_3Dnonlocal} \textbf{Non-locality of the weighted sum method for a dimension superior to two.}

Given the problem in (\ref{theo:moo:cex3D}), when the weighting vectors $\vec{w}^1 = (0.10, 0.10, 0.80)$, $\vec{w}^2 = (0.08, 0.85, 0.07)$ and $\vec{w}^3 = (0.32, 0.28, 0.40)$ are used in the weighted sum method, the following solutions are obtained: $\vec{r}^1 = (0.22, 0.22, 0.003)$, $\vec{r}^2 = (0.20, 0.001, 0.26)$, $\vec{r}^3 = (0.11, 0.15, 0.07)$. Then, the hyperplane with normal vector $\vec{w} = (0.16, 0.34, 0.50)$ is generated by a convex combination of the weighting vectors $\vec{w}^1$, $\vec{w}^2$ and $\vec{w}^3$. When this hyperplane is used in the weighted sum method, the solution $\vec{r} = (0.31209618,0.06910856,0.0318477)$ is found. It is possible to notice that the first component $r_1 = 0.31209618$ is larger than the maximum of the neighbors, $r^1_1 = 0.22$, demonstrating that the locality is not preserved.

\end{theorem}

\section{Properties of the weighting vector calculation in MONISE}\label{appendix:monise_ap}

In this appendix we prove that the problem of Definition \ref{def:moo:monise_w_estim} is bounded.

\begin{theorem} \label{theo:moo:monise_lim} The problem of Definition \ref{def:moo:monise_w_estim} is bounded.
\begin{proof} Let $(\vec{w}, \overline{\vec{r}}, \underline{\vec{r}})$ be any feasible solution of the problem. Then,

\begin{align}
& {\vec{w}}^\top\overline{\vec{r}} \leq {\vec{w}}^\top\vec{f}(\vec{x}^i),\\
& {\vec{w}}^\top\overline{\vec{r}}-{\vec{w}}^\top\underline{\vec{r}} \leq {\vec{w}}^\top\vec{f}(\vec{x}^i)-{\vec{w}}^\top\underline{\vec{r}} \leq {\vec{w}}^\top\vec{f}(\vec{x}^i) - {\vec{w}}^\top\vec{z}^{utopian}.
\end{align}
Since $\vec{z}^{utopian}$ and $\vec{f}(\vec{x}^i)$ are constants, the objective value ${\vec{w}}^\top\overline{\vec{r}}-{\vec{w}}^\top\underline{\vec{r}}$ we want to maximize is upper bounded. Therefore the problem is bounded.
\end{proof}

\end{theorem}

Then, we want to proof that, for bi-objective problems, the conditions imposed to the NISE neighborhood $(i,j)$ to determine $\vec{w}$ and $\underline{\vec{r}}$, related to that neighborhood, is such that the problem in Definition \ref{def:moo:monise_w_estim} satisfies the KKT conditions. First, we will prove some Lemmas to help in this demonstration.

\begin{lemma} \label{theo:moo:monise_aux1}
Given three distinct solutions $\vec{r}^i$, $\vec{r}^j$ and $\vec{r}^k \in \mathbb{R}^2$ obtained from the weighting vectors $\vec{w}^i$, $\vec{w}^j$ and $\vec{w}^k \in \mathbb{R}^2$, and assuming that $\vec{r}^k$ is not between $\vec{r}^i$ and $\vec{r}^j$, then there is no $\vec{p} \in \mathbb{R}^2$ such that:

\begin{equation} \label{eq:theo:monise_aux1:e1}
\begin{cases}
{\vec{w}^i}^\top \vec{p} - {\vec{w}^i}^\top\vec{r}^i = 0\\
{\vec{w}^j}^\top \vec{p} - {\vec{w}^j}^\top\vec{r}^j = 0\\
{\vec{w}^k}^\top \vec{p} - {\vec{w}^k}^\top\vec{r}^k = 0.
\end{cases}
\end{equation}

\begin{proof} Suppose, for example, that the third equality is true, and consider the following ordering (the opposite ordering will generate a similar proof): $w^i_1<w^j_1<w^k_1$ (consequently $w^i_2<w^j_2<w^k_2$).  By Theorem \ref{theo:moo:weight_2Dsens}, we have $r^i_1>r^j_1>r^k_1$ and $r^i_2<r^j_2<r^k_2$.

Rewrite $\vec{r}^j$ as $\vec{r}^j=\vec{p}+\alpha\begin{bmatrix} -w_2^j \\ w_1^j \end{bmatrix}$, $\alpha \in \mathbb{R}$, since ${\vec{w}^j}^\top \vec{r}^j={\vec{w}^j}^\top \vec{p}+\alpha\begin{bmatrix} w_1^j \\ w_2^j \end{bmatrix}^\top \begin{bmatrix} -w_2^j \\ w_1^j \end{bmatrix} = {\vec{w}^j}^\top \vec{p}$.\\

It is now possible to split it into two cases:

\begin{description}
\item[Case 1:] If $\alpha<0$, ${\vec{w}^i}^\top \vec{r}^j={\vec{w}^i}^\top \vec{p}+\alpha(-w_1^i w_2^j+w_2^i w_1^j)$, we have ${\vec{w}^i}^\top \vec{r}^j<{\vec{w}^i}^\top \vec{r}^i$ since $\alpha(-w_1^i w_2^j+w_2^i w_1^j)<0$.
\item[Case 2:] If $\alpha\geq 0$, ${\vec{w}^k}^\top \vec{r}^j={\vec{w}^k}^\top \vec{p}+\alpha(-w_1^k w_2^j+w_2^k w_1^j)$, we have ${\vec{w}^k}^\top \vec{r}^i\leq{\vec{w}^k}^\top \vec{r}^i$ since $\alpha(-w_1^k w_2^j+w_2^k w_1^j)<0$.
\end{description}
Since all cases generate a contradiction, the proof is concluded.
\end{proof}
\end{lemma}

\begin{lemma} \label{theo:moo:monise_aux2}
Given three solutions $\vec{r}^i$, $\vec{r}^j$ and $\vec{r}^k \in \mathbb{R}^2$, and assuming that $\vec{r}^k$ is not between $\vec{r}^i$ and $\vec{r}^j$, then there is no $\vec{p} \in \mathbb{R}^2$ such that:

\begin{equation} \label{eq:theo:monise_aux2:e2}
\begin{cases}
\vec{w}^\top \vec{p} - \vec{w}^\top\vec{r}^i = 0\\
\vec{w}^\top \vec{p} - \vec{w}^\top\vec{r}^j = 0\\
\vec{w}^\top \vec{p} - \vec{w}^\top\vec{r}^k = 0.
\end{cases}
\end{equation}

\begin{proof} Supposing by contradiction that the third equality is true.

Assuming the following ordering (the opposite ordering will generate a similar proof): $w^i_1<w^j_1<w^k_1$ (consequently $w^i_2<w^j_2<w^k_2$), using Theorem \ref{theo:moo:weight_2Dsens} we have that $r^i_1>r^j_1>r^k_1$ and $r^i_2<r^j_2<r^k_2$. First we can rewrite $\vec{r}^j$ in two ways $\vec{r}^j=\vec{r}^i+\alpha\begin{bmatrix} -w_2 \\ w_1 \end{bmatrix}, \alpha>0$ or $\vec{r}^j=\vec{r}^k+\beta\begin{bmatrix} w_2 \\ -w_1 \end{bmatrix}, \beta>0$. With the support of Theorem \ref{theo:moo:weight_2Drecur} we can see that $\vec{w}^\top \vec{p} - \vec{w}^\top\vec{r}^i = 0, \vec{w}^\top \vec{p} - \vec{w}^\top\vec{r}^k = 0$ will result in a weighting vector $\vec{w}$ between $\vec{w}^i$ and $\vec{w}^k$. Therefore we have two admissible situations:

\begin{description}
\item[Case 1:] $w_1^j<w_1$ will result in:
\begin{equation}
{\vec{w}^j}^\top\vec{r}^j={\vec{w}^j}^\top\vec{r}^k+\beta(w_2w_i^j-w_1w_2^j)>{\vec{w}^j}^\top\vec{r}^k.
\end{equation}
\item[Case 2:] $w_1^j>w_1$ will result in:
\begin{equation}
{\vec{w}^j}^\top\vec{r}^j={\vec{w}^j}^\top\vec{r}^i+\alpha(-w_2w_i^j+w_1w_2^j)>{\vec{w}^j}^\top\vec{r}^i.
\end{equation}
\end{description}
Since all cases generate a contradiction, the proof is concluded.
\end{proof}
\end{lemma}

\begin{theorem} \label{theo:moo:monise_nise} \textbf{NISE Equivalence}
The procedure adopted by NISE to find the weighting vector $\vec{w}$ and a reference point $\underline{\vec{r}}$ related to some neighborhood $(i,j)$ (Section \ref{sssec:moo:nise:choose}), satisfies the KKT conditions for the problem of Definition \ref{def:moo:monise_w_estim}.
\begin{proof}
The neighborhood $(i,j)$ is characterized by the solutions $\vec{r}^i$ and $\vec{r}^j$ found by the weighted sum method using weighting vectors $\vec{w}^i$ and $\vec{w}^j$. From Theorem \ref{theo:moo:weight_2Drecur}, $\vec{w} \geq \vec{0}$ and $\vec{w}^\top \vec{1} = 1$. Since $\vec{w}^i > 0$ and $\vec{w}^j > 0$ (the individual optimal solutions were found at initialization), it is easy to see that constraint  $\underline{\vec{r}} > \vec{z}^{utopian}$ is satisfied. Furthermore, as discussed in Section \ref{sssec:moo:nise:w_inf}, the constraints ${\vec{w}^i}^\top \vec{p} = {\vec{w}^i}^\top\vec{r}^i$ and ${\vec{w}^j}^\top \vec{p} = {\vec{w}^j}^\top\vec{r}^j$ are satisfied. By Lemma \ref{theo:moo:monise_aux1}, there is no other inequality of that type that is active. By what is demonstrated in Section \ref{sssec:moo:nise:w_inf}, the constraints $\vec{w}^\top \vec{p} = \vec{w}^\top\vec{r}^i$ and $\vec{w}^\top \vec{p} = \vec{w}^\top\vec{r}^j$ are satisfied. And due to Lemma \ref{theo:moo:monise_aux2} there is no other active inequality of that type.

Considering only active inequalities, the Lagrangian for the problem is written as:

\begin{multline}\label{eq:moo:lagrangean1}
l(\vec{w}, \overline{\vec{r}},\underline{\vec{r}}, \lambda, \kappa, \mu) = \vec{w}^\top \underline{\vec{r}} - \vec{w}^\top \overline{\vec{r}} + \lambda_i[{\vec{w}^i}^\top\vec{r}^i - {\vec{w}^i}^\top \underline{\vec{r}}] + \lambda_j[{\vec{w}^j}^\top\vec{r}^j - {\vec{w}^j}^\top \underline{\vec{r}}]+\\
+ \kappa_i [\vec{w}^\top \overline{\vec{r}} - \vec{w}^\top\vec{r}^i] + \kappa_j [\vec{w}^\top \overline{\vec{r}} - \vec{w}^\top\vec{r}^j] + \mu [\vec{w}^\top \vec{1} - 1].
\end{multline}
where $\lambda$, $\kappa$ and $\mu$ are Lagrange multipliers.

Applying the necessary conditions for optimality:

\begin{subequations}
\begin{align}
\nabla_\vec{w} l = \underline{\vec{r}} - \overline{\vec{r}} + \kappa_i [\overline{\vec{r}} - \vec{r}^i] + \kappa_j [\overline{\vec{r}} - \vec{r}^j] + \mu \vec{1} = \vec{0},\\
\nabla_{\overline{\vec{r}}} l = - \vec{w} +  \kappa_i \vec{w} + \kappa_j \vec{w} = \vec{0},\\
\nabla_{\underline{\vec{r}}} l = \vec{w} - \lambda_i \vec{w}^i - \lambda_j \vec{w}^j = \vec{0}.
\end{align}
\end{subequations}
we conclude that $\lambda_i, \lambda_j, \kappa_i, \kappa_j > 0$. Then we can see that the NISE method satisfies the necessary conditions for the proposed problem.
\end{proof}
\end{theorem}

It is worth noting that the objective function $\vec{w}^\top \overline{\vec{r}}-\vec{w}^\top\underline{\vec{r}}$ is different from the scalar $\mu^{i,j} = \frac{[\vec{w}^\top \vec{r}^i-\vec{w}^\top\underline{\vec{r}}]}{||\vec{w}||}$ which guides the choice of neighborhood in NISE. The recursive steps (sequence of neighborhoods) in MONISE may be different from that of NISE. The objective function adopted by MONISE simplifies the optimization model, leading to an equivalent mixed-integer linear formulation.

\section{A mixed-integer linear equivalent model for weighting vector calculation in MONISE} \label{appendix:monise:lip}

We solve the non-convex problem in Definition \ref{def:moo:monise_w_estim}, by solving an equivalent mixed-integer linear problem, which allows the use of commercial or open-source solvers. To do so, we employed a procedure, inspired by \cite{Xia2015}, that uses KKT conditions to find an equivalent mixed-integer linear problem.

Replacing $\vec{w}^\top \overline{\vec{r}}$ with $v$ and writing the Lagrangian, we have:

\begin{multline}\label{eq:moo:lagrangean2}
l(\vec{w}, \overline{\vec{r}},\underline{\vec{r}}, v, \lambda, \kappa, \gamma, \beta, \nu, \mu) = \vec{w}^\top \underline{\vec{r}} - v - \sum_{i=1}^L \lambda_i[{\vec{w}^i}^\top \underline{\vec{r}} - {\vec{w}^i}^\top\vec{r}^i] +\\
+ \sum_{i=1}^L \kappa_i [v - \vec{w}^\top\vec{r}^i] + \sum_{i=1}^m \beta_i [\underline{\vec{r}}_i-\vec{z}^{utopian}_i] -\sum_{i=1}^m \nu_i\vec{w}_i + \mu [\vec{w}^\top \vec{1} - 1].
\end{multline}

Aiming at creating a new equivalent formulation, some KKT conditions are applied:

\begin{equation} \label{eq:moo:KKTw}
\nabla_\vec{w} l = \underline{\vec{r}} - \sum_{i=1}^L \kappa_i \vec{r}^i  - \nu + \mu \vec{1} = 0,
\end{equation}

\begin{equation} \label{eq:moo:KKTv}
\nabla_v l = -1 + \sum_{i=1}^L \kappa_i = 0,
\end{equation}

\begin{equation} \label{eq:moo:KKTineq}
\kappa_i [-(v - \vec{w}^\top\vec{r}^i)] = 0\ \forall i \in \{1,\ldots,L\},
\end{equation}

\begin{equation} \label{eq:moo:KKTwineq}
\nu_i w_i = 0 \ \forall i \in \{1,\ldots,m\}.
\end{equation}

Taking the inner product of both sides of (\ref{eq:moo:KKTw}) with $\vec{w}$, we have:

\begin{equation} \label{eq:moo:KKTinnerw}
\vec{w}^\top\underline{\vec{r}} - \sum_{i=1}^L \kappa_i \vec{w}^\top\vec{r}^i  - \vec{w}^\top\nu + \mu \vec{w}^\top\vec{1} = 0.
\end{equation}

From (\ref{eq:moo:KKTv}) we have that $\sum_{i=1}^L \kappa_i v - v = 0$. Summing up $\sum_{i=1}^L \kappa_i v - v$ in (\ref{eq:moo:KKTinnerw}) we obtain:

\begin{equation} \label{eq:moo:KKTinnerw2}
\vec{w}^\top\underline{\vec{r}} + \sum_{i=1}^L \kappa_i\left(v - \vec{w}^\top\vec{r}^i\right) -v  - \vec{w}^\top\nu + \mu \vec{w}^\top\vec{1} = 0.
\end{equation}

Using KKT conditions and the original constraints, (\ref{eq:moo:KKTinnerw2}) produces:

\begin{equation} \label{eq:moo:fobjnew}
\vec{w}^\top\underline{\vec{r}} -v = -\mu.
\end{equation}

Since the complementary KKT conditions in (\ref{eq:moo:KKTineq}) are true when $\kappa_i = 0$ or $[-(v - \vec{w}^\top\vec{r}^i)]=0$, and (\ref{eq:moo:KKTwineq}) is true when $\nu_i = 0$ or $w_i = 0$, it is possible to express the same conditions of (\ref{eq:moo:KKTineq}) using linear-integer constrains such as:

\begin{subequations}\label{eq:moo:KKTineqLI}
\begin{align}
[-(v - \vec{w}^\top\vec{r}^i)] \geq 0,\\
\kappa_i \geq 0,\\
[-(v - \vec{w}^\top\vec{r}^i)] \leq \kappa_i^B \overline{\kappa_i},\\
\kappa_i \leq (1-\kappa_i^B),
\end{align}
\end{subequations}
where $\kappa_i^B$ is a binary variable that is equal to one when $\kappa_i=0$ and equal to zero when $[-(v - \vec{w}^\top\vec{r}^i)]=0$, and $\overline{\kappa_i}$ is a constant always larger than $[-(v - \vec{w}^\top\vec{f}(\vec{x}^{i})]$. Similar constraints are used to express (\ref{eq:moo:KKTwineq}):

\begin{subequations}\label{eq:moo:KKTwineqLI}
\begin{align}
w_i \geq 0,\\
\nu_i \geq 0,\\
w_i \leq \nu_i^B,\\
\nu_i \leq (1-\nu_i^B)\overline{\nu_i},
\end{align}
\end{subequations}
where $\nu_i^B$ is a binary variable that is equal to one when $\nu_i=0$ and equal to zero when $w_i=0$ and $\overline{\nu_i}$ is a constant always larger than $\nu_i$.

Changing the objective function to $\mu$, expressed by (\ref{eq:moo:fobjnew}), and adding KKT conditions in (\ref{eq:moo:KKTw}) and (\ref{eq:moo:KKTv}), as well as the linear--integer equivalents in (\ref{eq:moo:KKTineqLI}) and (\ref{eq:moo:KKTwineqLI}), we reach an equivalent mixed-integer linear problem for the problem in Definition \ref{def:moo:monise_w_estim}, which is formally stated in Definition \ref{def:moo:monise_w_estim_pli}.
\end{appendices}
\end{document}